\documentclass[a4paper,11pt]{article}
\usepackage[utf8]{inputenc}
\usepackage{fullpage}
\usepackage{amsmath}
\usepackage{amsfonts}
\usepackage{mathtools}
\usepackage{bbm}
\usepackage{amssymb}
\usepackage{graphicx}
\usepackage{epstopdf}
\usepackage{caption}
\usepackage{color}
\usepackage{subcaption}
\usepackage{siunitx}
\usepackage{setspace}
\usepackage{psfrag}
\usepackage[square, comma, numbers, compress, sort]{natbib}
\usepackage{rotating}
\usepackage{authblk}
\usepackage{multirow}
\usepackage[displaymath, mathlines]{lineno}

\newcommand*\patchAmsMathEnvironmentForLineno[1]{%
  \expandafter\let\csname old#1\expandafter\endcsname\csname #1\endcsname
  \expandafter\let\csname oldend#1\expandafter\endcsname\csname end#1\endcsname
  \renewenvironment{#1}%
     {\linenomath\csname old#1\endcsname}%
     {\csname oldend#1\endcsname\endlinenomath}}%
\newcommand*\patchBothAmsMathEnvironmentsForLineno[1]{%
  \patchAmsMathEnvironmentForLineno{#1}%
  \patchAmsMathEnvironmentForLineno{#1*}}%
\AtBeginDocument{%
\patchBothAmsMathEnvironmentsForLineno{equation}%
\patchBothAmsMathEnvironmentsForLineno{align}%
\patchBothAmsMathEnvironmentsForLineno{flalign}%
\patchBothAmsMathEnvironmentsForLineno{alignat}%
\patchBothAmsMathEnvironmentsForLineno{gather}%
\patchBothAmsMathEnvironmentsForLineno{multline}%
}

\title{A power series method for solving ordinary and partial differentials equations motivated by domain growth.}
\author[1]{Robert J. H. Ross \thanks{robert\_ross@hms.harvard.edu}}
\affil[1]{Harvard Medical School, Systems Biology Department, Warren Alpert Building, 25 Shattuck St, Boston, MA 02115}
\begin{document}

\maketitle

\begin{abstract}
\noindent In this work we present a power series method for solving ordinary and partial differential equations. To demonstrate our method we solve a system of ordinary differential equations describing the movement of a random walker on a one-dimensional lattice, two nonlinear ordinary differential equations, a wave and diffusion equation (linear partial differential equations), and a nonlinear partial differential equation (quasilinear). The inclusion of boundary conditions and the general solutions to other equations of interest are included in the Supplementary material.
\end{abstract}

{\begin{center}{\noindent{\bf Keywords:}
Differential equations, ordinary, partial, power series solutions, growth.}
\end{center}}

\section{Introduction}

We present a method that can be used to generate power series solutions for linear and nonlinear ordinary differential equations (ODEs) and partial differential equations (PDEs). The method we present relies on a separation of variables in a system of equations we construct, and generates a power series weighted by coefficients written in terms of the initial condition.  This method was motivated by work conducted on the mathematical modeling of domain growth \citep{Ross2017a}. 
\\
\\
The outline of this work is as follows: In Section 2.1 we demonstrate how our method can be used to generate power series solutions for linear and nonlinear ODEs.  In Sections 2.2 and 2.3 we demonstrate how this method can be extended to linear and nonlinear PDEs, including the implementation of boundary conditions in the linear PDE case.  We finish with a short discussion of the method presented in this work in Section 3.

\section{Results}

\subsection{Solving a system of ordinary differential equations describing the movement of a random walker on a one-dimensional lattice}

\begin{figure}[h!]
\centering
\includegraphics[width=0.3\textwidth]{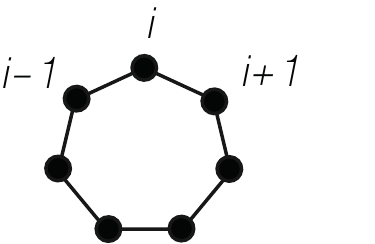}
\caption{A one-dimensional lattice with periodic boundary conditions can be represented as a ring.  The sites are sequentially labelled from $i \in \{1, 2, ..., N \}$, with $N$ being the total number of sites.}
\label{fig:figure1}
\end{figure}

\begin{figure}[h!]
\centering
\begin{subfigure}[b]{0.4\textwidth}
	\includegraphics[width=1\textwidth]{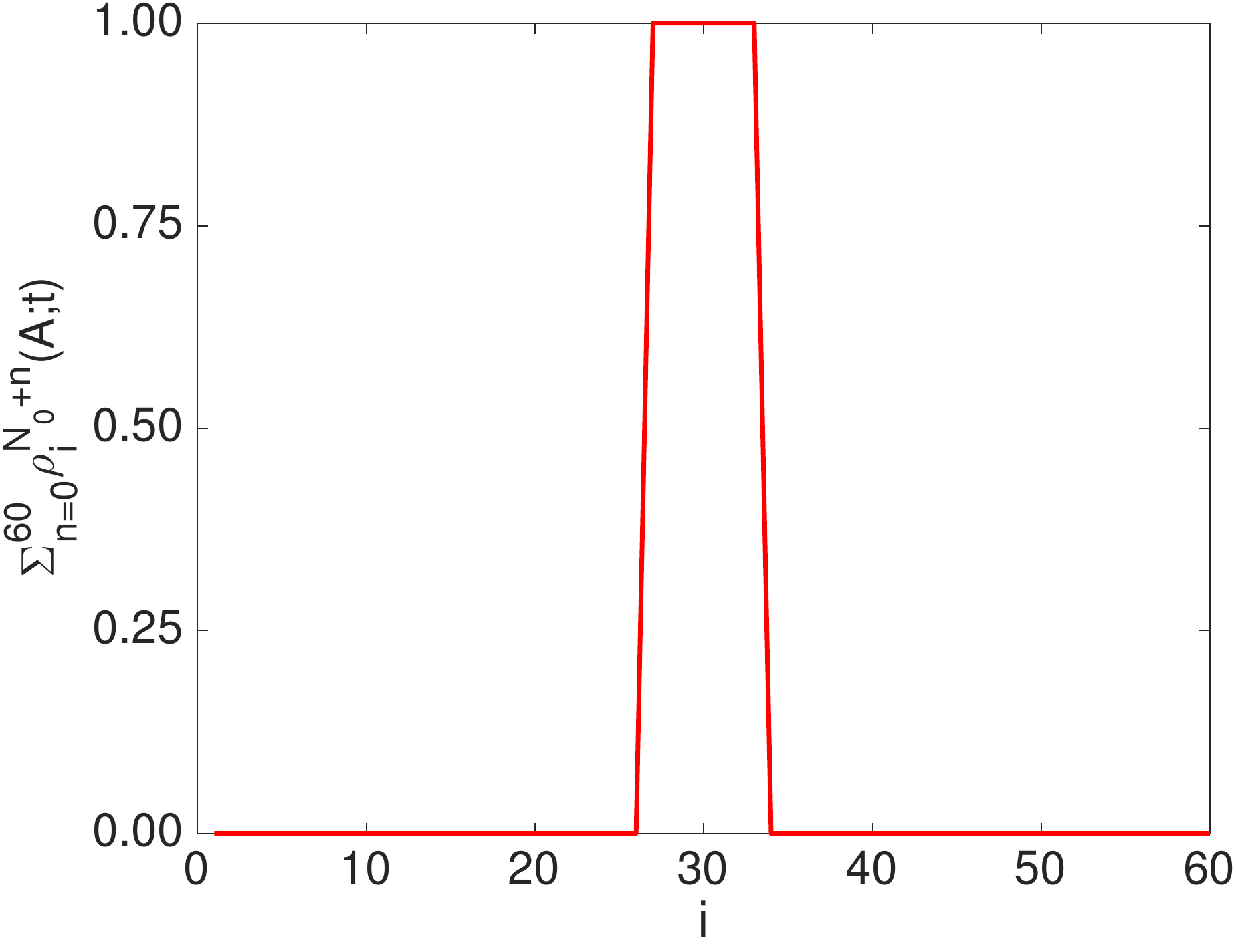}
	\subcaption{}
\end{subfigure}
\begin{subfigure}[b]{0.4\textwidth}
	\includegraphics[width=1\textwidth]{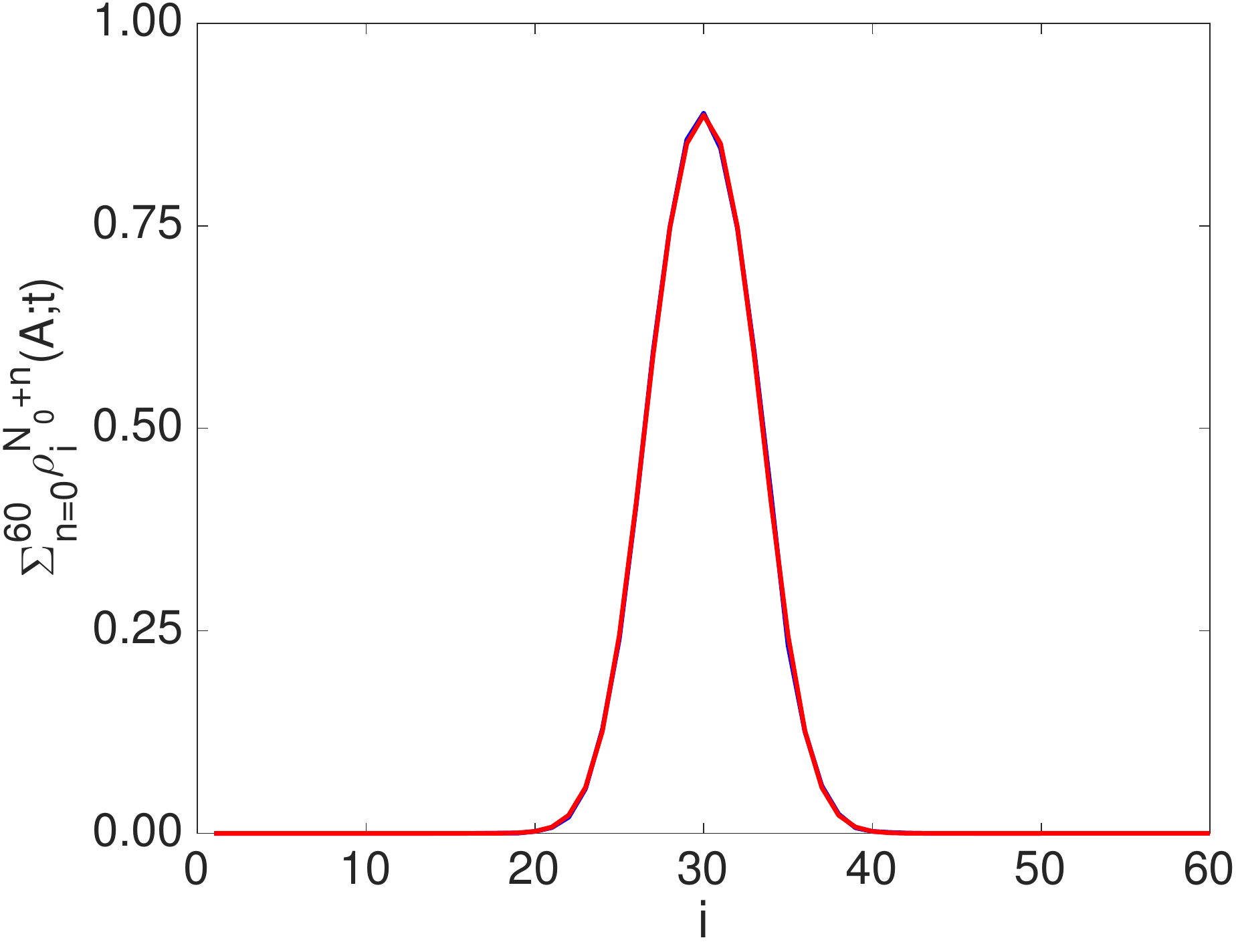}
	\subcaption{}
\end{subfigure}
\begin{subfigure}[b]{0.4\textwidth}
	\includegraphics[width=1\textwidth]{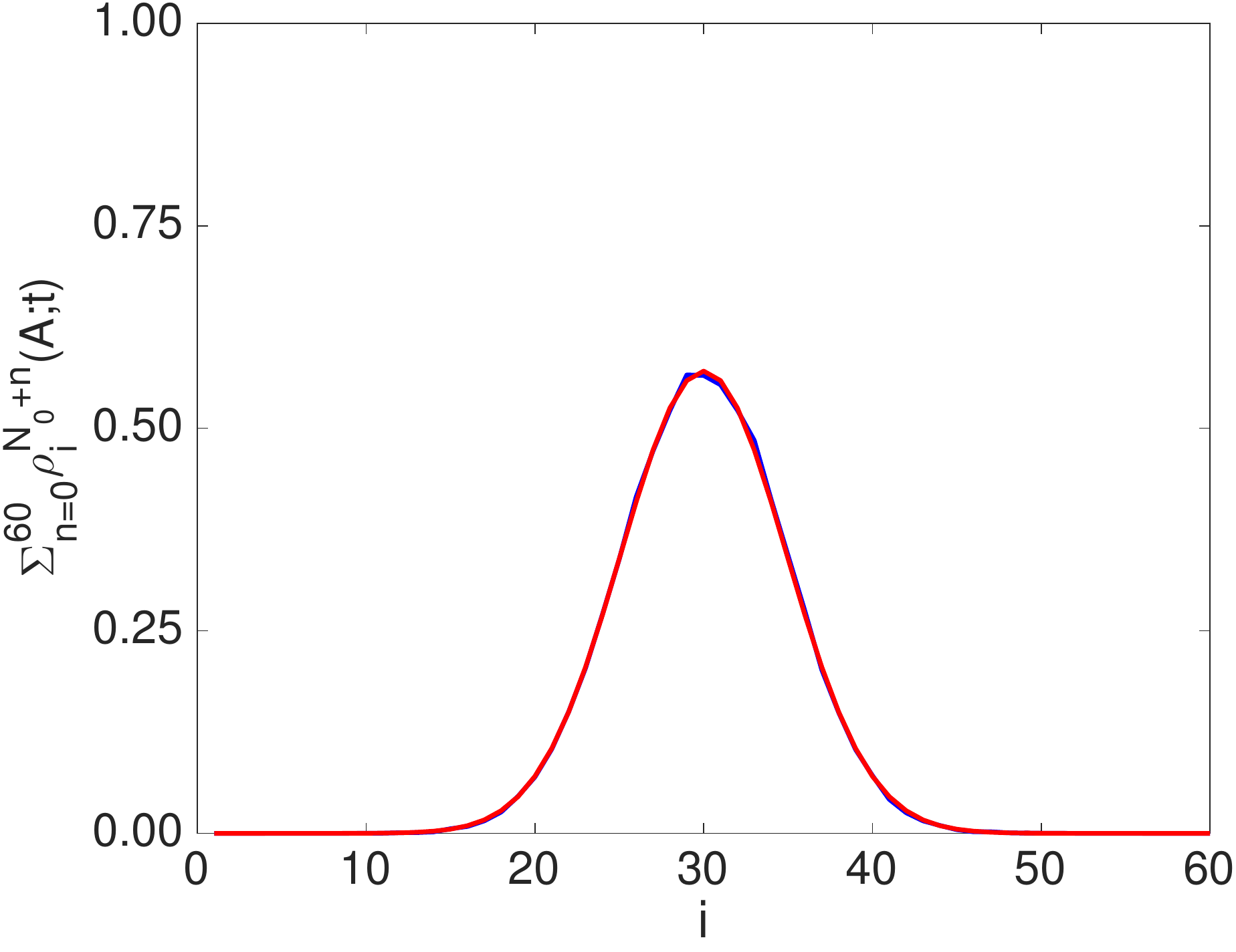}
	\subcaption{}
\end{subfigure}
\caption{A comparison of an ensemble average of the discrete model with periodic boundary conditions and Eq. \eqref{eq:metastatic_diffusion_sum} at different time points. The blue lines indicate the ensemble average from the discrete model and the red lines indicate the solutions of Eq. \eqref{eq:metastatic_diffusion_sum}. In the discrete model agents were placed from sites 27:33 for the initial condition in each replicate. The ensemble average was calculated from 10000 replicates of the discrete model, and Eq. \eqref{eq:metastatic_diffusion_sum} was truncated at $n=60$.  In the discrete model and Eq. \eqref{eq:metastatic_diffusion_sum} $P_{m} = 1$ and $\beta = 1$. In (a) $t = 0$, in (b) $t = 50$, and in (c) $t = 200$.}
\label{fig:figure2}
\end{figure}

A one-dimensional lattice with periodic boundaries (a ring) is displayed in Fig. \ref{fig:figure1}. 
\\
The following equation describes the time evolution of the probability that an unbiased excluding random walker occupies site $i$ on a one-dimensional periodic lattice\footnote{A derivation of Eq. \eqref{eq:static_diffusion} can be found in the Supplementary material (SM1).}:
\begin{align}
\frac{\mathrm{d}p_{i}}{\mathrm{d}t} = \left(\frac{P_{m}}{2}\right)\left(p_{i-1} - 2p_{i} + p_{i+1}\right).
\label{eq:static_diffusion}
\end{align}
In Eq. \eqref{eq:static_diffusion} $p_{i}$ is the probability a random walker is situated at site $i$ at time $t$, and $P_{m}$ is the rate at which the random walker attempts to move to an adjacent site on the lattice.
\\
\\
We now write Eq. \eqref{eq:static_diffusion} in the following manner
\begin{align}
\frac{\mathrm{d}}{\mathrm{d}t}\sum^{\infty}_{n=0}p^{n}_{i} = \sum^{\infty}_{n=0}\left(\frac{P_{m}}{2}\right)\left(p^{n}_{i-1} - 2p^{n}_{i} + p^{n}_{i+1}\right).
\label{eq:sum_back}
\end{align}
That is, we postulate that $p_{i}$ can be written as the infinite series
\begin{align} 
p_{i} = \sum^{\infty}_{n=0}p^{n}_{i}.
\label{eq:pbar}
\end{align}
We now decompose Eq. \eqref{eq:sum_back} into the following infinite system of equations
\begin{align}
\frac{\mathrm{d}p^{0}_{i}}{\mathrm{d}t} = -\beta p^{0}_{i},
\label{eq:metastatic_diffusion_asep_0}
\end{align}
and
\begin{align}
\frac{\mathrm{d}p^{n}_{i}}{\mathrm{d}t} = -\beta p^{n}_{i} + \beta p^{n-1}_{i} + \left(\frac{P_{m}}{2}\right)\left(p^{n-1}_{i-1} - 2p^{n-1}_{i} + p^{n-1}_{i+1}\right), \ \ \ \forall \ n > 0.
\label{eq:metastatic_diffusion_asep}
\end{align}
Notice in Eq. \eqref{eq:metastatic_diffusion_asep} that we have written the terms associated with the movement of the random walker in terms of $n-1$, not $n$.\footnote{It has previously been shown \citep{Ross2017a} that if the lattice is growing, Eq. \eqref{eq:metastatic_diffusion_asep} would be written as
\begin{align}
\frac{\mathrm{d}p^{n}_{i}}{\mathrm{d}t} = -P_{g} p^{n}_{i} + P_{g} p^{n-1}_{i} + \left(\frac{P_{m}}{2}\right)\Big(p^{n}_{i-1} - 2p^{n}_{i} + p^{n}_{i+1}\Big), \nonumber
\end{align}
where the motility terms have the same $n$ as the time derivative, and $P_{g}$ is the rate at which the lattice grows.  This observation is what initially motivated the work we present here.}
We include the parameter $\beta$ in Eqs. \eqref{eq:metastatic_diffusion_asep_0} and \eqref{eq:metastatic_diffusion_asep} as a `shape' parameter, however, its inclusion in Eqs. \eqref{eq:metastatic_diffusion_asep_0} and \eqref{eq:metastatic_diffusion_asep} is not necessary and $\beta$ can be set to zero if desired.\footnote{If $\beta = 0$ then Eqs. \eqref{eq:metastatic_diffusion_asep_0} and \eqref{eq:metastatic_diffusion_asep} would be
\begin{align}
\frac{\mathrm{d}p^{0}_{i}}{\mathrm{d}t} = 0,
\end{align}
and
\begin{align}
\frac{\mathrm{d}p^{n}_{i}}{\mathrm{d}t} = \left(\frac{P_{m}}{2}\right)\left(p^{n-1}_{i-1} - 2p^{n-1}_{i} + p^{n-1}_{i+1}\right), \ \ \ \forall \ n > 0.
\end{align}} Finally, we simplify Eq. \eqref{eq:metastatic_diffusion_asep} to obtain
\begin{align}
\frac{\mathrm{d}p^{n}_{i}}{\mathrm{d}t} = -\beta p^{n}_{i} + (\beta - P_{m})p^{n-1}_{i} + \left(\frac{P_{m}}{2}\right)\left(p^{n-1}_{i-1} + p^{n-1}_{i+1}\right), \ \ \ \forall \ n > 0.
\label{eq:metastatic_diffusion2}
\end{align}
It is readily apparent that if we sum Eq. \eqref{eq:metastatic_diffusion_asep_0} and Eq. \eqref{eq:metastatic_diffusion_asep} for all $n > 0$ we obtain
\begin{align}
\frac{\mathrm{d}}{\mathrm{d}t}\sum^{\infty}_{n=0}p^{n}_{i} = \sum^{\infty}_{n=0}\left(\frac{P_{m}}{2}\right)\Big(p^{n}_{i-1} - 2p^{n}_{i} + p^{n}_{i+1}\Big).
\label{eq:metastatic_diffusion_typ2}
\end{align}
This means if we substitute Eq. \eqref{eq:pbar} into Eq. \eqref{eq:metastatic_diffusion_typ2} we arrive at
\begin{align}
\frac{\mathrm{d} p_{i}}{\mathrm{d}t} = \left(\frac{P_{m}}{2}\right)\left(p_{i-1} - 2p_{i} + p_{i+1}\right),
\label{eq:sum_back_sub}
\end{align}
which recapitulates Eq. \eqref{eq:static_diffusion}.
\\
\\
The decomposition of Eq. \eqref{eq:static_diffusion} into the infinite system of equations contained in Eqs. \eqref{eq:metastatic_diffusion_asep_0} and \eqref{eq:metastatic_diffusion_asep} is straightforward to solve.  To see this consider the initial equation, Eq. \eqref{eq:metastatic_diffusion_asep_0},
\begin{align}
\frac{\mathrm{d}p^{0}_{i}}{\mathrm{d}t} = -\beta p^{0}_{i}.
\label{eq:metastatic_diffusion_init}
\end{align}
As the initial equation does not `inherit' any terms associated with the movement of the random walker, Eq. \eqref{eq:metastatic_diffusion_init} admits the simple solution:
\begin{align}
p^{0}_{i} = A_{i}\operatorname{exp}(-\beta t),
\label{eq:metastatic_diffusion_init_sol}
\end{align}
where $A_{i}$ is the initial value of site $i$.
\\
\\
Equation \eqref{eq:metastatic_diffusion_init_sol} can be placed in Eq. \eqref{eq:metastatic_diffusion2} for $n = 1$, so that Eq. \eqref{eq:metastatic_diffusion2} becomes
\begin{align}
\frac{\mathrm{d}p^{1}_{i}}{\mathrm{d}t} = -\beta p^{1}_{i} + (\beta - P_{m})A_{i}\operatorname{exp}(-\beta t) + \left(\frac{P_{m}}{2}\right)\left(A_{i-1} + A_{i+1}\right)\operatorname{exp}(-\beta t).
\label{eq:metastatic_diffusion_simp}
\end{align}
This means Eq. \eqref{eq:metastatic_diffusion_simp} is now also straightforward to solve. Repeated application of this process admits the following recurrence formula as a solution for Eq. \eqref{eq:metastatic_diffusion2}:
\begin{align}
p_{i}^{n}(t) = \operatorname{exp}(-\beta t)\left(\frac{t^{n}}{n!}\right)\left[ \sum_{j=0}^{n}(\beta -P_{m})^{n-j}(P_{m})^{j}\frac{1}{2^{j}}\binom{n}{j}\left(\sum_{k=-j}^{-j:2:j}\binom{j}{\frac{(k+j)}{2}}A_{i+k}\right)\right], \ \ \ \forall \ n \geq 0.
\label{eq:metastatic_diffusion_sol}
\end{align}
Therefore, the probability of site $i$ being occupied by a walker at time $t$ is given by
\begin{align}
\sum_{n = 0}^{\infty} p_{i}^{n}(t) = \sum_{n = 0}^{\infty}\operatorname{exp}(-\beta t)\left(\frac{t^{n}}{n!}\right)\left[ \sum_{j=0}^{n}(\beta -P_{m})^{n-j}(P_{m})^{j}\frac{1}{2^{j}}\binom{n}{j}\left(\sum_{k=-j}^{-j:2:j}\binom{j}{\frac{(k+j)}{2}}A_{i+k}\right)\right],
\label{eq:metastatic_diffusion_sum}
\end{align}
in accordance with Eq. \eqref{eq:pbar}.
\\
\\
In Fig. \ref{fig:figure2} we demonstrate that the solution given by Eq. \eqref{eq:metastatic_diffusion_sum} matches the evolution of the ensemble average of the discrete model excellently at our chosen level of truncation (for the algorithm used in the discrete model see the Supplementary material (SM2)).  In Fig. \ref{fig:figure2b} we display what we refer to as `streams', given by Eq. \eqref{eq:metastatic_diffusion_sol}, for a single value of the shape parameter $\beta $.  In the Supplementary material (SM3) we display streams for different values of $\beta$, which demonstrates how $\beta $ influences the shape of the streams that compose the solution given by Eq. \eqref{eq:metastatic_diffusion_sum}.

\begin{figure}[h!]
\centering
\begin{subfigure}[b]{0.4\textwidth}
	\includegraphics[width=1\textwidth]{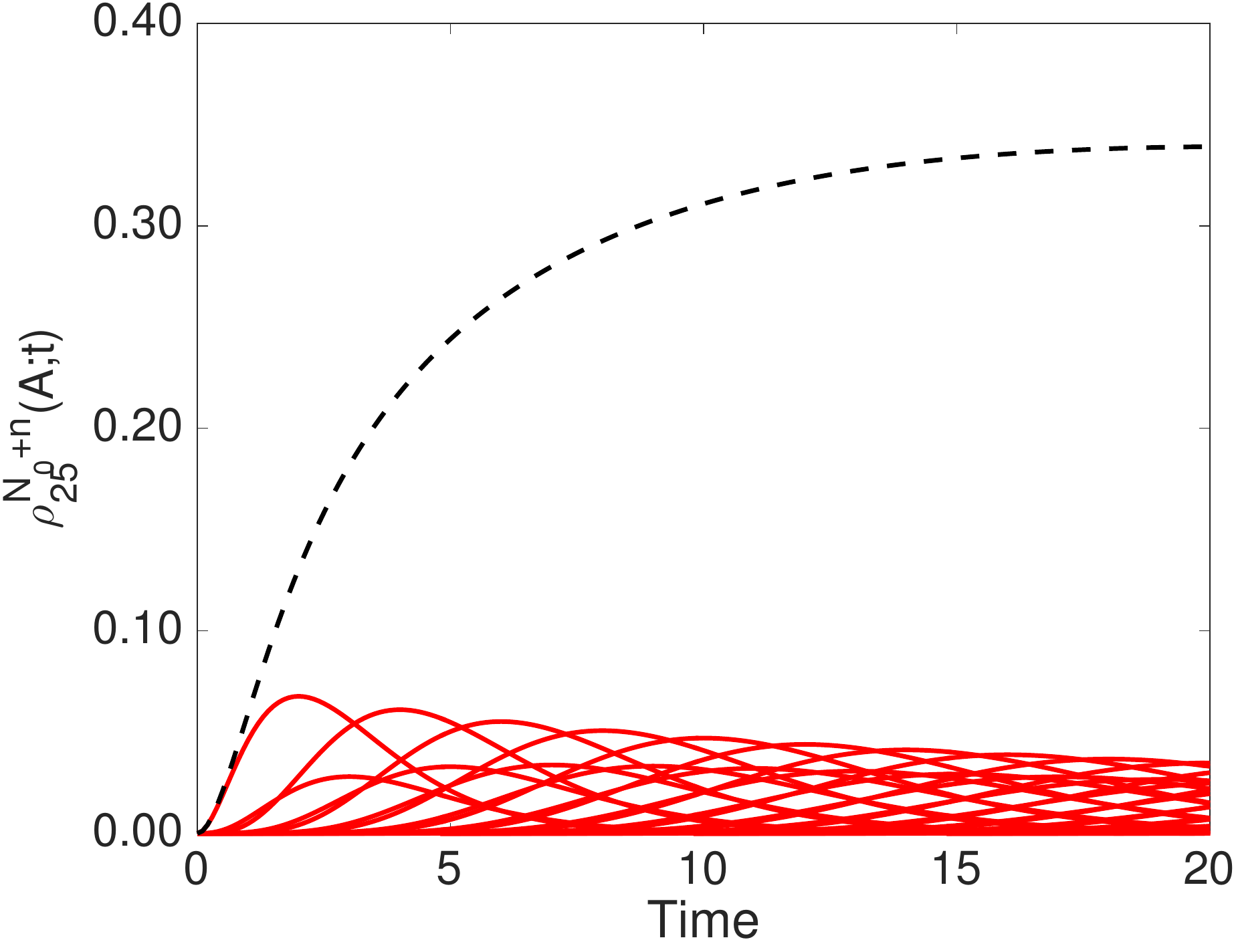}
	\subcaption{}
\end{subfigure}
\begin{subfigure}[b]{0.4\textwidth}
	\includegraphics[width=1\textwidth]{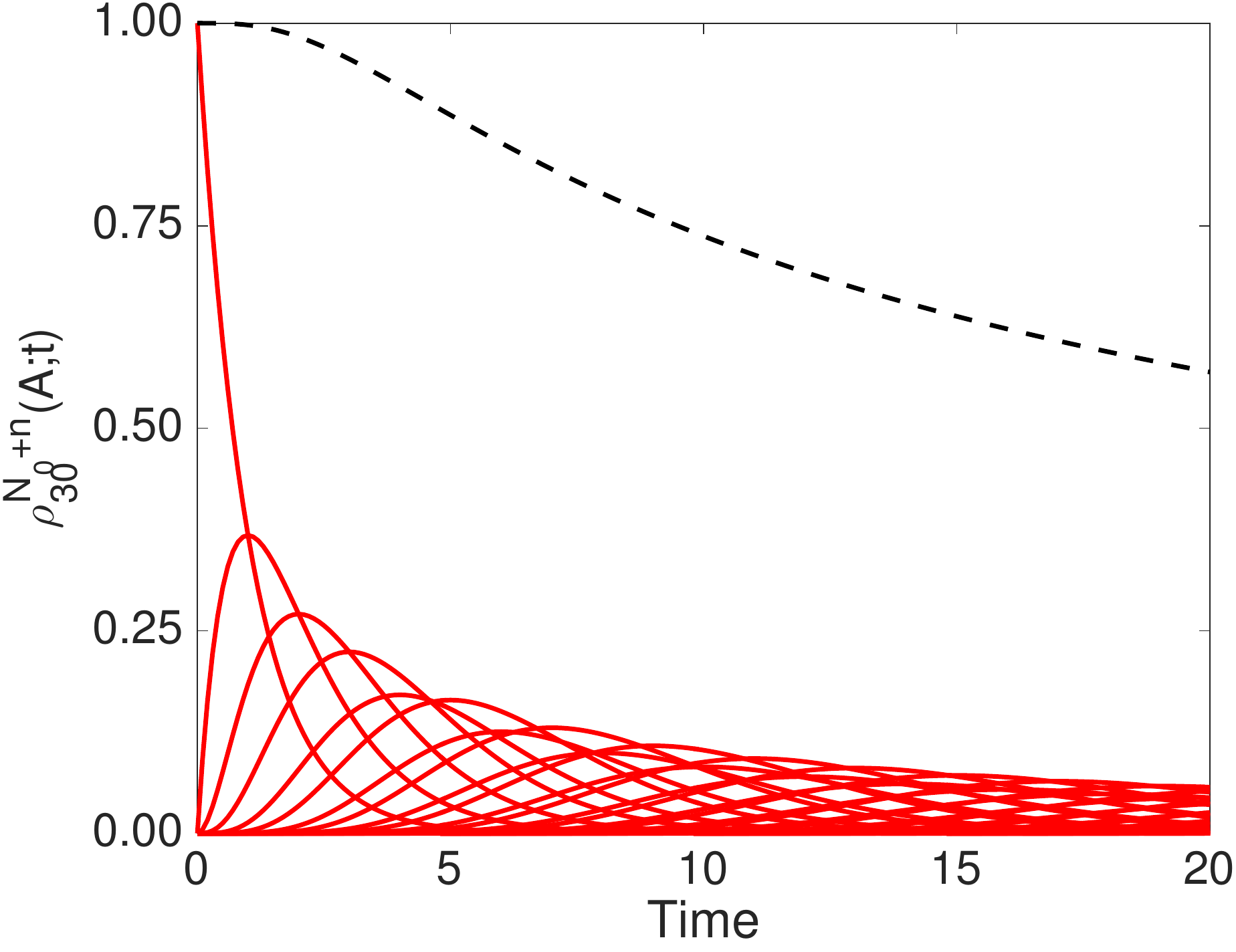}
	\subcaption{}
\end{subfigure}
\caption{The streams of sites $i = 25$ and $i=30$ as given by Eq. \eqref{eq:metastatic_diffusion_sol} for incrementing values of $n$ from 1:60, where $n$ increases from the left to the right in both panels (a) and (b). For both panels (a) and (b) $P_{m} = 1$ and $\beta = 1$. In (a) $i = 25$, and in (b) $i = 30$. The black-dashed line is the sum of all the streams for site $i$, given by Eq. \eqref{eq:metastatic_diffusion_sum}. }
\label{fig:figure2b}
\end{figure}

\subsubsection{Boundary conditions}

\begin{figure}[h!]
\centering
\begin{subfigure}[b]{0.4\textwidth}
	\includegraphics[width=1\textwidth]{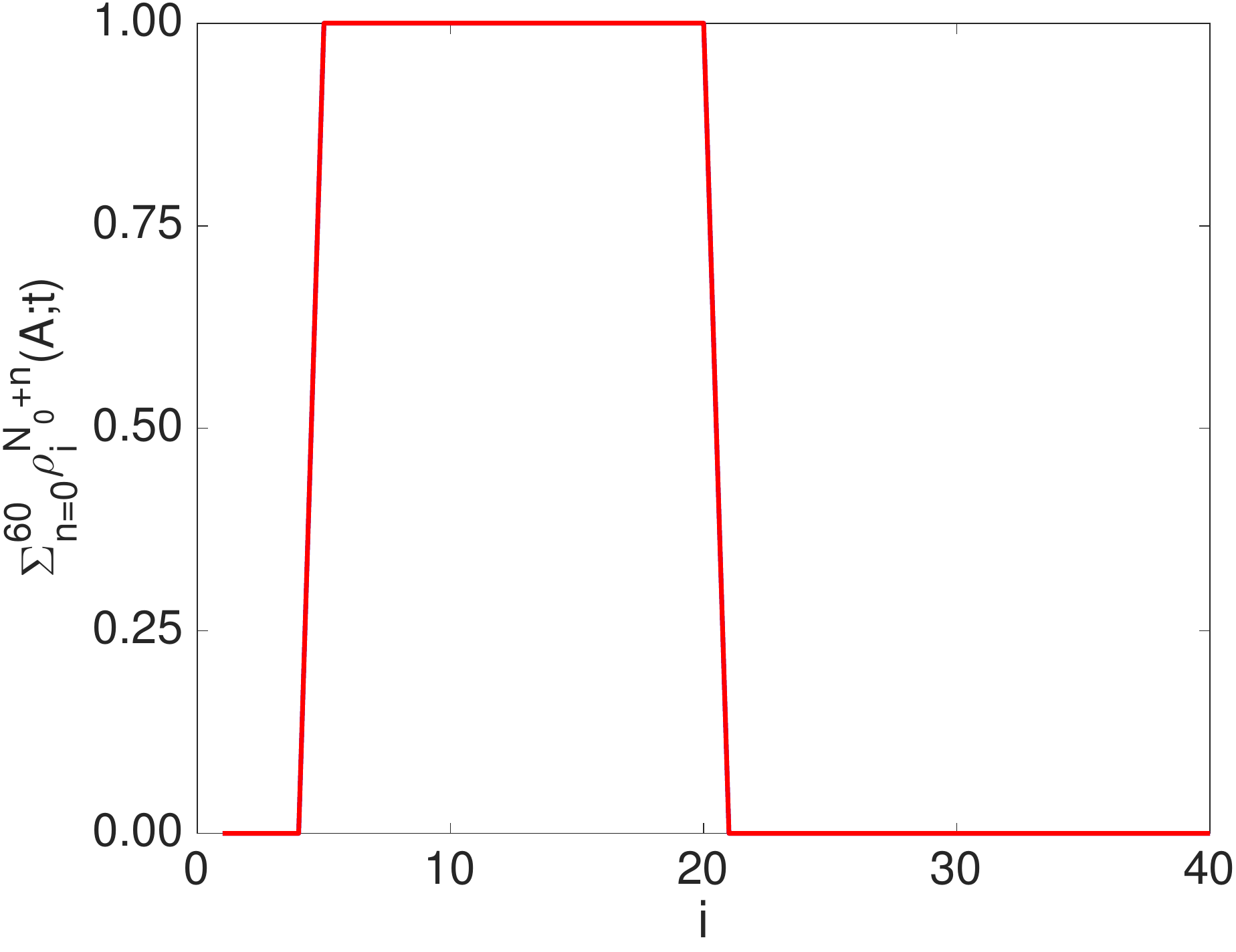}
	\subcaption{}
\end{subfigure}
\begin{subfigure}[b]{0.4\textwidth}
	\includegraphics[width=1\textwidth]{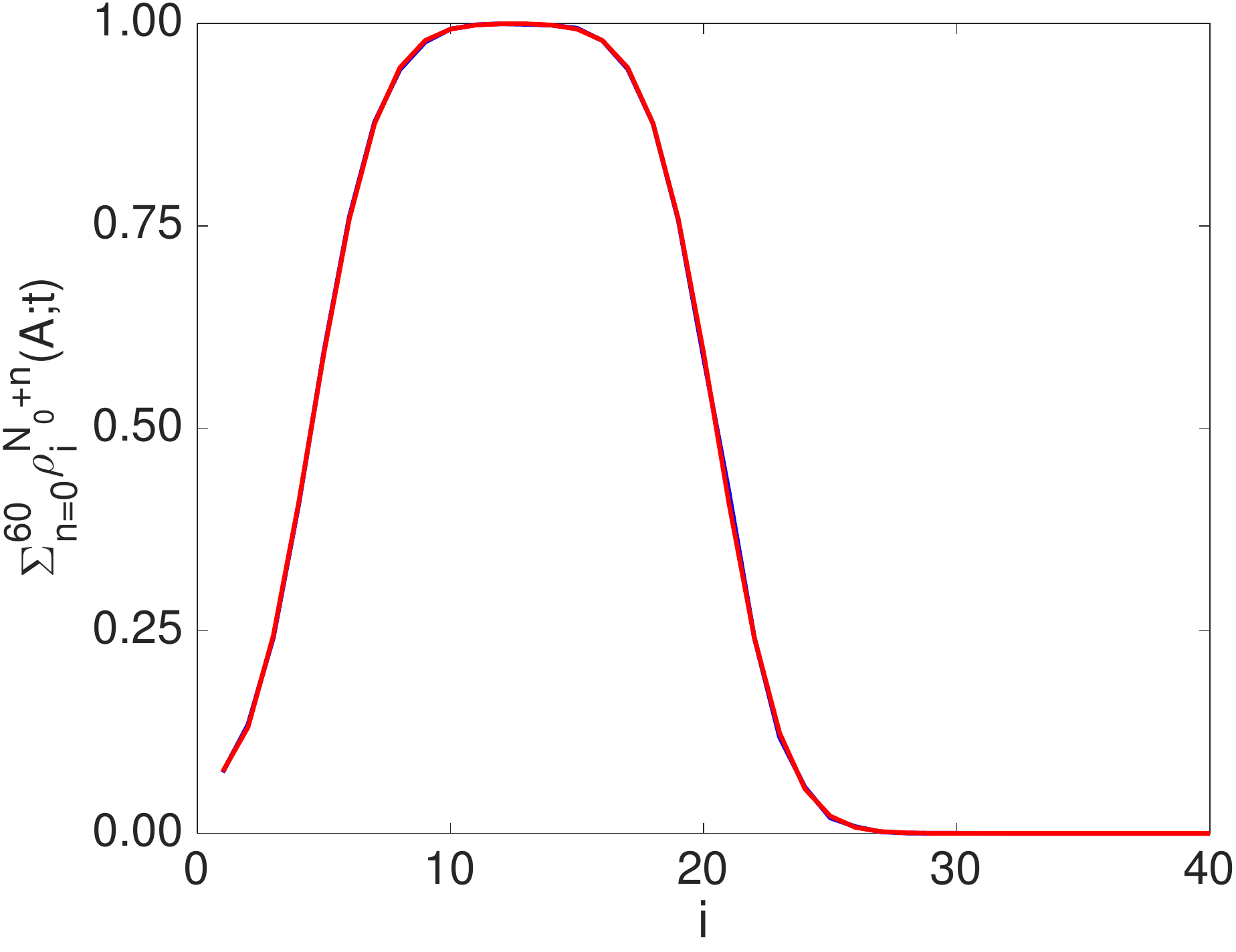}
	\subcaption{}
\end{subfigure}
\begin{subfigure}[b]{0.4\textwidth}
	\includegraphics[width=1\textwidth]{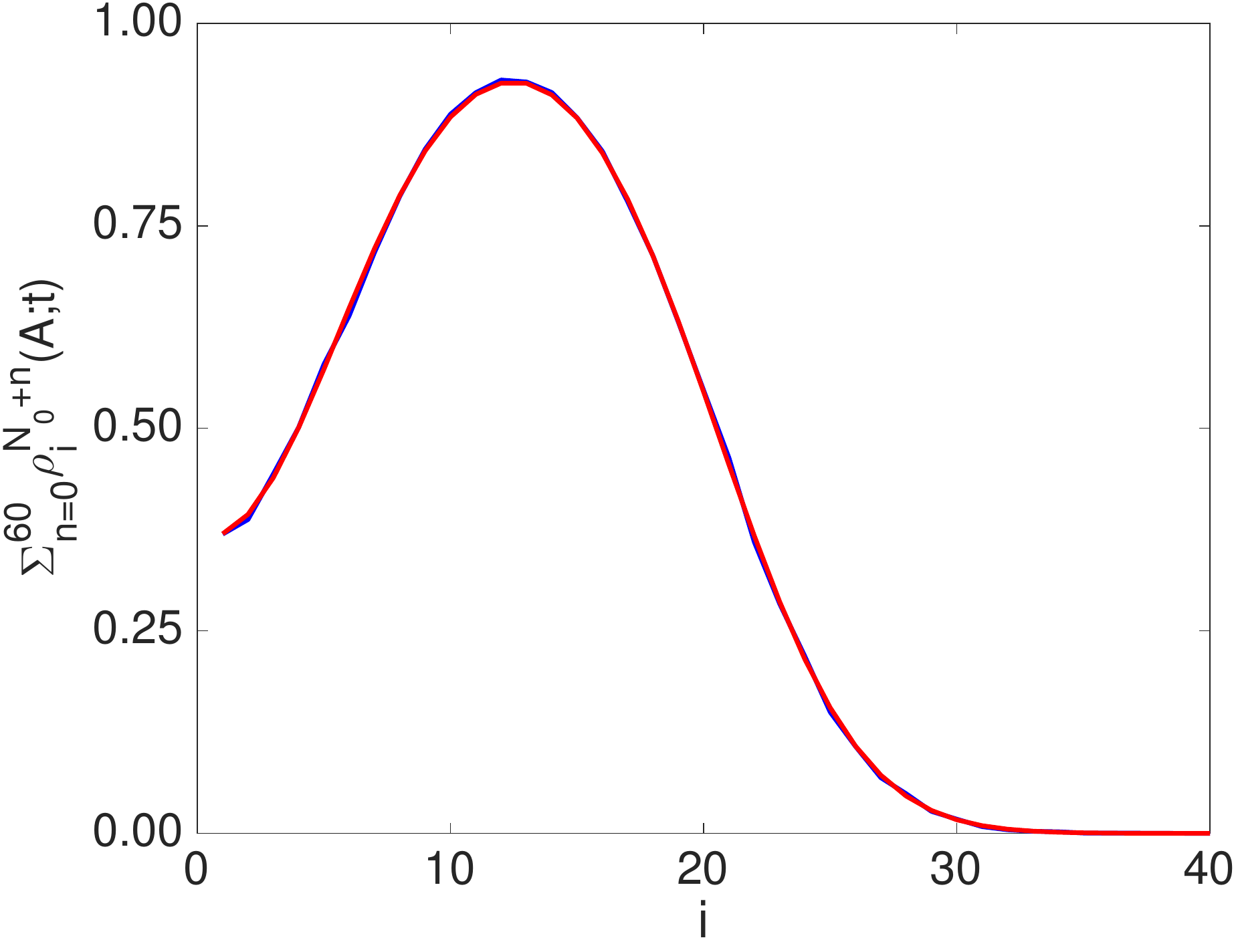}
	\subcaption{}
\end{subfigure}
\caption{A comparison of an ensemble average of the discrete model with no-flux boundary conditions and Eqs. \eqref{eq:middle}, \eqref{eq:bound1gen} and \eqref{eq:bound2gen} at different time points. The blue lines indicate the ensemble average and the red lines indicate Eqs. \eqref{eq:middle}, \eqref{eq:bound1gen} and \eqref{eq:bound2gen}. In the discrete model agents were placed from sites 5:20 for the initial condition in each replicate. The ensemble average was calculated from 10000 replicates of the discrete model, and Eqs. \eqref{eq:middle}, \eqref{eq:bound1gen} and \eqref{eq:bound2gen} are truncated at $n = 100$.  In the discrete model and Eq. \eqref{eq:metastatic_diffusion_sum} $P_{m} = 1$. In (a) $t = 0$, in (b) $t = 50$, and in (c) $t = 200$.}
\label{fig:figure3}
\end{figure}

It is possible to implement boundary conditions in equations describing the movement of a random walker on a one-dimensional lattice with the method we are presenting. For an unbiased random walker on a one-dimensional lattice with no-flux boundary conditions the equations describing the probability of finding a walker at a given site are:
\begin{align}
\frac{\mathrm{d}p^{n}_{1}}{\mathrm{d}t} &= -\beta p^{n}_{1} + \beta p^{n-1}_{1} + \left(\frac{P_{m}}{2}\right)\left(p^{n-1}_{2} - p^{n-1}_{1}\right), \nonumber \\ 
&\vdotswithin{=} \nonumber \\
\frac{\mathrm{d}p^{n}_{i}}{\mathrm{d}t} &= -\beta p^{n}_{i} + \beta p^{n-1}_{i} + \left(\frac{P_{m}}{2}\right)\left(p^{n-1}_{i-1} - 2p^{n-1}_{i} + p^{n-1}_{i+1}\right), \nonumber \\
&\vdotswithin{=} \nonumber \\
\frac{\mathrm{d}p^{n}_{N}}{\mathrm{d}t} &= -\beta p^{n}_{N} + \beta p^{n-1}_{N} + \left(\frac{P_{m}}{2}\right)\left(p^{n-1}_{N-1} - p^{n-1}_{N}\right), \ \ \ \forall \ n > 0.
\label{eq:metastatic_diffusionBC1}
\end{align}
Following the same procedure as we did for Eqs. \eqref{eq:metastatic_diffusion_sol} and \eqref{eq:metastatic_diffusion_sum} we find the following recurrence relation for sites not situated on the boundary
\begin{align}
p_{i}^{n}(t) = \left(\frac{t}{n}\right)\left[(\beta -P_{m})p_{i}^{n-1}(t) + \left(\frac{P_{m}}{2}\right)(p_{i+1}^{n-1}(t) + p_{i-1}^{n-1}(t))\right], \ \ \ \forall \ n > 0,
\label{eq:ecom}
\end{align}
with 
\begin{align}
p_{i}^{0}(t) = A_{i}\operatorname{exp}(-\beta t).
\label{}
\end{align}
The reader will notice that we have written Eq. \eqref{eq:ecom} in a more economical form than Eq. \eqref{eq:metastatic_diffusion_sol}.  It is possible to write Eq. \eqref{eq:ecom} in the same manner as Eq. \eqref{eq:metastatic_diffusion_sol}, and if done so each site on the lattice will have its own recurrence formula describing the probability of a random walker being located at that site at time $t$.
\\
\\
From Eq. \eqref{eq:ecom} the evolution with respect to time of the probability of site $i$ being occupied by a walker is given by
\begin{align}
\sum^{\infty}_{n=0}p_{i}^{n}(t) = p_{i}^{0}(t) + \sum^{\infty}_{n=1}\left(\frac{t}{n}\right)\left[(\beta -P_{m})p_{i}^{n-1}(t) + \left(\frac{P_{m}}{2}\right)(p_{i+1}^{n-1}(t) + p_{i-1}^{n-1}(t))\right].
\label{eq:middle}
\end{align}
The recurrence relations for sites situated on the boundary are
\begin{align}
p_{1}^{n}(t) = \left(\frac{t}{n}\right)\left[\left(\beta -\frac{P_{m}}{2}\right)p_{1}^{n-1}(t) + \left(\frac{P_{m}}{2}\right)(p_{2}^{n-1}(t))\right] \ \ \ \forall \ n > 0,
\label{eq:bound1}
\end{align}
and
\begin{align}
p_{N}^{n}(t) = \left(\frac{t}{n}\right)\left[\left(\beta -\frac{P_{m}}{2}\right)p_{N}^{n-1}(t) + \left(\frac{P_{m}}{2}\right)(p_{N-1}^{n-1}(t))\right] \ \ \ \forall \ n > 0.
\label{eq:bound2}
\end{align}
Therefore, the evolution with respect to time of the probability of sites $1$ and $N$ being occupied by a walker are given by
\begin{align}
\sum^{\infty}_{n=0}p_{1}^{n}(t) = p_{1}^{0}(t) + \sum^{\infty}_{n=1}\left(\frac{t}{n}\right)\left[\left(\beta -\frac{P_{m}}{2}\right)p_{1}^{n-1}(t) + \left(\frac{P_{m}}{2}\right)(p_{2}^{n-1}(t))\right],
\label{eq:bound1gen}
\end{align}
and
\begin{align}
\sum^{\infty}_{n=0}p_{N}^{n}(t) = p_{N}^{0}(t) + \sum^{\infty}_{n=1}\left(\frac{t}{n}\right)\left[\left(\beta -\frac{P_{m}}{2}\right)p_{N}^{n-1}(t) + \left(\frac{P_{m}}{2}\right)(p_{N-1}^{n-1}(t))\right],
\label{eq:bound2gen}
\end{align}
respectively.
\\
\\
In Fig. \ref{fig:figure3} the solutions given by Eqs. \eqref{eq:middle}, \eqref{eq:bound1gen} and \eqref{eq:bound2gen} are displayed.  It can be seen that Eqs. \eqref{eq:middle}, \eqref{eq:bound1gen} and \eqref{eq:bound2gen} and the ensemble average from the discrete model match excellently.  The algorithm for the discrete model can be found in the Supplementary material (SM2).
\\
\\
We provide a final example of our method being applied to a linear ODE in the Supplementary material (SM4).

\subsubsection{Nonlinear ordinary differential equations}

\begin{figure}[h!]
\centering
\begin{subfigure}[b]{0.4\textwidth}
	\includegraphics[width=1\textwidth]{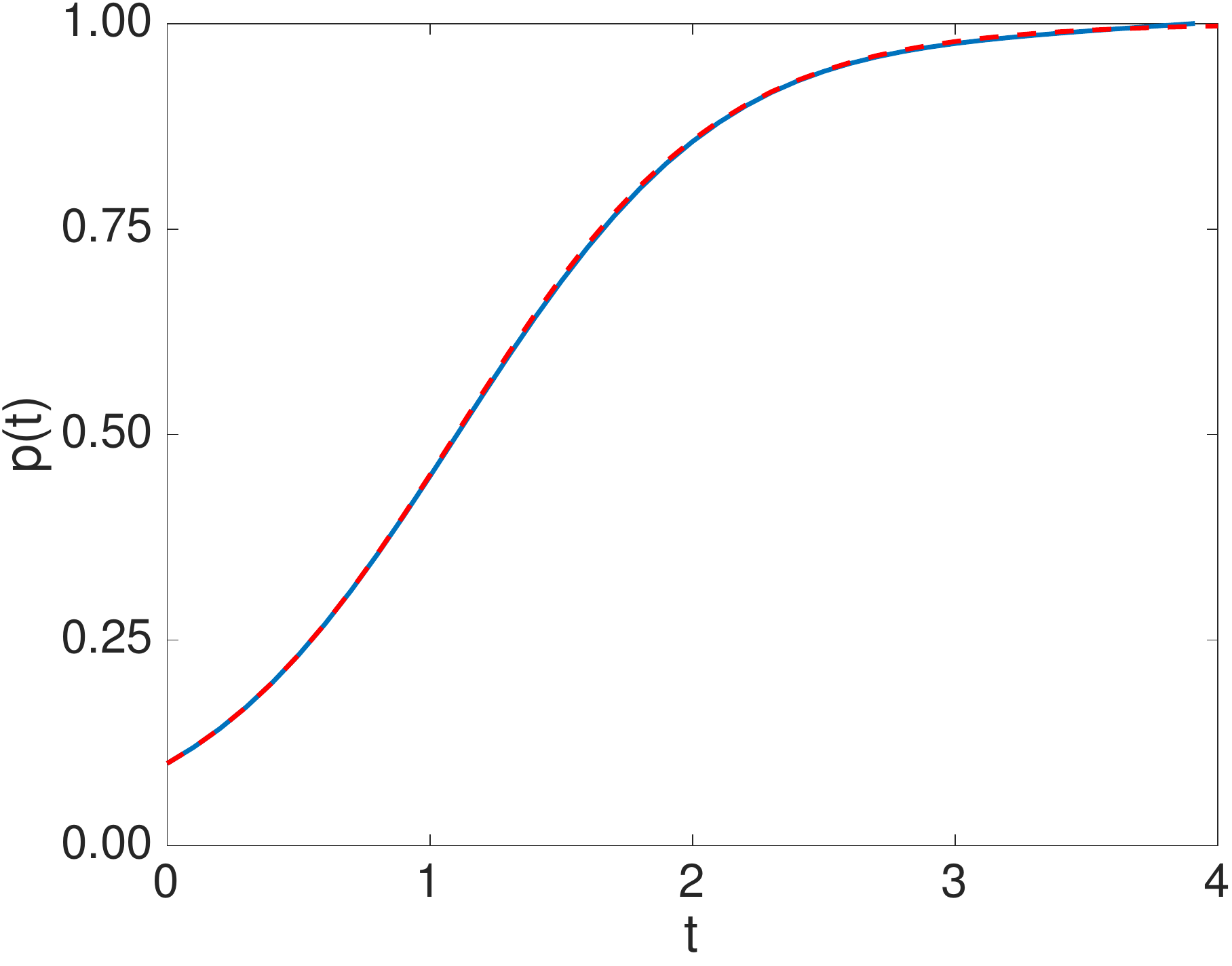}
	\subcaption{}
\end{subfigure}
\begin{subfigure}[b]{0.4\textwidth}
	\includegraphics[width=1\textwidth]{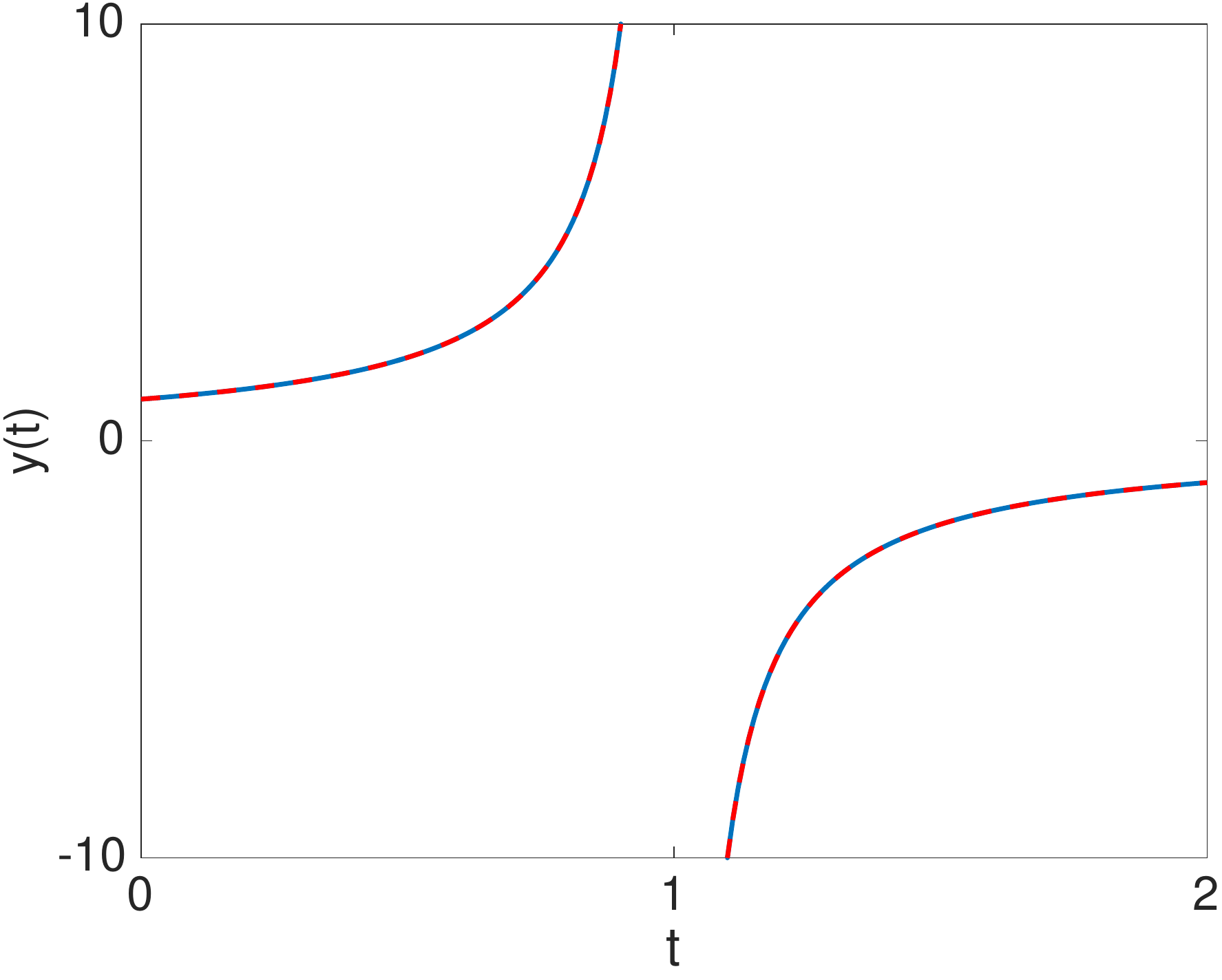}
	\subcaption{}
\end{subfigure}
\caption{In (a) Eqs. \eqref{eq:nonlin_exp_N0_new1} and \eqref{eq:nonlin_exp_N0_new4} are compared with Eq. \eqref{eq:non_lin2_sol} for $\gamma = 3$ and $C_{1} = 0.1$. The truncation value for Eq. \eqref{eq:nonlin_exp_N0_new4} in (a) is $n = 10$ and $\beta = 10$.  We use symbolic integration in MATLAB to solve Eqs. \eqref{eq:nonlin_exp_N0_new1} and \eqref{eq:nonlin_exp_N0_new4}. In (b) Eq. \eqref{eq:non_lin3sol} is compared with Eq. \eqref{eq:non_lin3anal} for $y^{0} = 1$ and $\alpha = 1$. The truncation value for Eq. \eqref{eq:non_lin3sol} in (b) is $n=20$.}
\label{fig:figureNLode}
\end{figure}

We now apply our method to nonlinear ODEs.  This allows us to demonstrate another important aspect of our method.
Initially we solve
\begin{align}
\frac{\mathrm{d}p}{\mathrm{d}t} = \gamma p(1-p),
\label{eq:non_lin2}
\end{align}
where $\gamma > 0$. The analytic solution of Eq. \eqref{eq:non_lin2} is
\begin{align}
p(t) = \frac{C_{1}\operatorname{exp}(\gamma t)}{1 - C_{1} + C_{1}\operatorname{exp}(\gamma t)},
\label{eq:non_lin2_sol}
\end{align}
where $C_{1}$ is the value of $p(t)$ at $t = 0$.\\
\\
To solve Eq. \eqref{eq:non_lin2} in our framework we proceed as follows. To begin with we decompose Eq. \eqref{eq:non_lin2} into
\begin{align}
\frac{\mathrm{d}p^{n}(t)}{\mathrm{d}t} = -\beta p^{n}(t) + \beta p^{n-1}(t) + \gamma p^{n-1}(t) - \gamma p^{n-1}(t)\sum^{\infty}_{n=0}p^{n}(t), \ \ \ \forall \ n > 0,
\label{eq:nonlin_exp3}
\end{align}
with
\begin{align}
\frac{\mathrm{d}p^{0}(t)}{\mathrm{d}t} = -\beta p^{0}(t).
\label{eq:nonlin_exp4}
\end{align}
It is evident that Eq. \eqref{eq:nonlin_exp3} cannot be solved in the same iterative manner as Eq. \eqref{eq:metastatic_diffusion_asep} due to the due to the nonlinear term present on its right-hand-side.\footnote{One might think the Eq. \eqref{eq:nonlin_exp3} should be written as
\begin{align}
\frac{\mathrm{d}p^{n}(t)}{\mathrm{d}t} = -\beta p^{n}(t) + \beta p^{n-1}(t) + \gamma p^{n-1}(t) - \gamma p^{n-1}(t)p^{n-1}(t), \ \ \ \forall \ n > 0,
\label{eq:nonlin_exp_wrong}
\end{align}
but this is incorrect as each stream needs to be multiplied by every other stream to account for the nonlinearity in Eq. \eqref{eq:non_lin2}.} To circumvent this we sum Eq. \eqref{eq:nonlin_exp3} for all $n \geq 0$ to obtain
\begin{align}
\sum^{\infty}_{n=0}\frac{\mathrm{d}p^{n}(t)}{\mathrm{d}t} = \gamma\sum^{\infty}_{n=0}p^{n}(t) - \gamma\sum^{\infty}_{n=0}p^{n}(t)\sum^{\infty}_{n=0}p^{n}(t),
\label{eq:expand}
\end{align}
and then decompose Eq. \eqref{eq:expand} in the following manner:
\begin{align}
\frac{\mathrm{d}p^{0}(t)}{\mathrm{d}t} = -\beta p^{0}(t),
\label{eq:nonlin_exp_N0_new1}
\end{align}
\begin{align}
\frac{\mathrm{d}p^{1}(t)}{\mathrm{d}t} = -\beta p^{1}(t) + \beta p^{0}(t) + \gamma p^{0} - \gamma p^{0}(t)p^{0}(t),
\label{eq:nonlin_exp_N0_new2}
\end{align}
\begin{align}
\frac{\mathrm{d}p^{2}(t)}{\mathrm{d}t} = -\beta p^{2}(t) + \beta p^{1}(t) + \gamma p^{1} - \gamma p^{1}(t)p^{1}(t) - 2\gamma p^{1}(t)p^{0}(t),
\label{eq:nonlin_exp_N0_new3}
\end{align}
so that in general
\begin{align}
\frac{\mathrm{d}p^{n}(t)}{\mathrm{d}t} &= -\beta p^{n}(t) + \beta p^{n-1}(t) + \gamma p^{n-1}(t) \nonumber \\
& \ \ \ - \gamma p^{n-1}(t)p^{n-1}(t) - 2\gamma p^{n-1}(t)\left(\sum^{n-2}_{i=0}p^{i}(t)\right), \ \ \ \forall n > 0.
\label{eq:nonlin_exp_N0_new4}
\end{align}
The decomposition of Eq. \eqref{eq:expand} into the equations contained in Eqs. \eqref{eq:nonlin_exp_N0_new1} and \eqref{eq:nonlin_exp_N0_new4} allows us to solve the unknowns iteratively, and it is straightforward to demonstrate that summing Eqs. \eqref{eq:nonlin_exp_N0_new1} and \eqref{eq:nonlin_exp_N0_new4} for all $n > 0$ returns Eq. \eqref{eq:expand}.  In Fig. \ref{fig:figureNLode} (a) we compare the solution of Eq. \eqref{eq:nonlin_exp_N0_new4} with the analytical solution Eq. \eqref{eq:non_lin2_sol}.
\\
It is also possible to derive power series solutions in terms of simple functions for nonlinear ODEs with the method we are presenting. For instance, the nonlinear ODE
\begin{align}
\frac{\mathrm{d}y}{\mathrm{d}t} = \alpha y^{2},
\label{eq:non_lin3}
\end{align}
has the following power series solution
\begin{align}
y = \sum^{\infty}_{n=0} \frac{(-1)^{n}y^{0}}{n!}\operatorname{log}^{n}(1 - \alpha y^{0} t),
\label{eq:non_lin3sol}
\end{align}
where
\begin{align}
y^{0} = y(0) = A.
\label{eq:non_lin3yo}
\end{align}
In Fig. \ref{fig:figureNLode} (b) we compare Eq. \eqref{eq:non_lin3sol} with the analytical solution of Eq. \eqref{eq:non_lin3},
\begin{align}
y = \frac{1}{\frac{1}{A} - \alpha t}.
\label{eq:non_lin3anal}
\end{align}
The details of how to derive Eq. \eqref{eq:non_lin3sol} are given in the Supplementary material (SM5). 

\subsection{Solving a linear partial differential equation}

We now demonstrate that the method we are presenting is also applicable to PDEs.
We start by applying this method to linear PDEs, for instance the wave equation:
\begin{align}
\frac{\partial u(x,t)}{\partial t} = c\frac{\partial u(x,t)}{\partial x}.
\label{eq:hyperbolic_eq}
\end{align}
Motivated by the previous section we write Eq. \eqref{eq:hyperbolic_eq} as
\begin{align}
\frac{\partial u^{n}}{\partial t} = -\beta u^{n} + \beta u^{n-1} + c\frac{\partial u^{n-1}}{\partial x}, \ \ \ \forall n > 0,
\label{eq:hyperbolic_eq_growth}
\end{align}
with
\begin{align}
\frac{\partial u^{0}}{\partial t} = -\beta u^{0}.
\label{eq:hyperbolic_eq_growth0}
\end{align}
As before we initially solve Eq. \eqref{eq:hyperbolic_eq_growth0},
\begin{align}
u^{0}(x,t) = A(x)\operatorname{exp}(-\beta t).
\label{eq:pde_diffusion_init}
\end{align}
It can be seen in Eq. \eqref{eq:pde_diffusion_init} that our method relies on the separation of spatial and temporal variable in the initial equation.
For general $u^{n}$ we obtain
\begin{align}
u^{n}(x,t) = \left(\frac{t^{n}}{n!}\right)\operatorname{exp}(-\beta t)\left[\sum^{n}_{j=0}\binom{n}{j}(\beta )^{n-j}c^{j}\left(\frac{\partial^{j}A(x)}{\partial x^{j}} \right)\right], \ \ \ \forall \ n \geq 0.
\label{eq:hyperbolic_eq_sol}
\end{align}
If $\beta  = 0$ we have
\begin{align}
u^{n}(x,t) = \left(\frac{t^{n}}{n!}\right)\left[c^{n}\left(\frac{\partial^{n}A(x)}{\partial x^{n}} \right)\right], \ \ \ \forall \ n \geq 0,
\label{eq:hyperbolic_eq_sol_Pg0}
\end{align}
and each stream is a polynomial in $t$ weighted by coefficients written in terms of the initial condition\footnote{In the case of Eq. \eqref{eq:hyperbolic_eq_sol_Pg0} it is evident we have simply derived a Taylor series expansion \citep{Abbott2001}.}.
From Eq. \eqref{eq:hyperbolic_eq_sol} the general solution to Eq. \eqref{eq:hyperbolic_eq} is
\begin{align}
\sum^{\infty}_{n=0}u^{n}(x,t) = \sum^{\infty}_{n=0}\left(\frac{t^{n}}{n!}\right)\operatorname{exp}(-\beta t)\left[\sum^{n}_{j=0}\binom{n}{j}(\beta )^{n-j}c^{j}\left(\frac{\partial^{j}A(x)}{\partial x^{j}} \right)\right].
\label{eq:hyperbolic_eq_sol_gen}
\end{align}
A simple initial condition for Eq. \eqref{eq:hyperbolic_eq_sol_gen} is
\begin{align}
A(x) = \operatorname{sin}(x),
\label{eq:ic}
\end{align}
and if we substitute Eq. \eqref{eq:ic} into Eq. \eqref{eq:hyperbolic_eq_sol} we obtain
\begin{align}
\sum^{\infty}_{n=0}u^{n}(x,t) = \sum^{\infty}_{n=0}\left(\frac{t^{n}}{n!}\right)\operatorname{exp}(-\beta t)\left[\sum^{n}_{j=0}\binom{n}{j}(\beta )^{n-j}c^{j}\left(\frac{\partial^{j}}{\partial x^{j}}\operatorname{sin}{x} \right)\right].
\label{eq:hyperbolic_eq_sol_sin}
\end{align}
It is also straightforward to solve the diffusion equation, which is
\begin{align}
\frac{\partial u}{\partial t} = D\frac{\partial^{2} u}{\partial x^{2}}.
\label{eq:diffusion_eq}
\end{align}
Solving Eq. \eqref{eq:diffusion_eq} in a similar manner to how we solved Eq. \eqref{eq:hyperbolic_eq} we obtain
\begin{align}
u^{n}(x,t) = \left(\frac{t^{n}}{n!}\right)\operatorname{exp}(-\beta t)\left[\sum^{n}_{j=0}\binom{n}{j}(\beta )^{n-j}D^{j}\left(\frac{\partial^{2j}A(x)}{\partial x^{2j}}\right)\right], \ \ \ \forall \ n \geq 0.
\label{eq:diffusion_eq_sol}
\end{align}
If $\beta  = 0$ Eq. \eqref{eq:diffusion_eq_sol} is
\begin{align}
u^{n}(x,t) = \left(\frac{t^{n}}{n!}\right)\left[D^{n}\left(\frac{\partial^{2n}A(x)}{\partial x^{2n}}\right)\right], \ \ \ \forall \ n \geq 0.
\label{eq:diffusion_eq_sol_Pg0}
\end{align}
Therefore, from Eq. \eqref{eq:diffusion_eq_sol} our general solution to Eq. \eqref{eq:diffusion_eq} is
\begin{align}
\sum_{n=0}^{\infty}u^{n}(x,t) &= \sum_{n=0}^{\infty}\left(\frac{t^{n}}{n!}\right)\operatorname{exp}(-\beta t)\left[\sum^{n}_{j=0}\binom{n}{j}(\beta )^{n-j}D^{j}\left(\frac{\partial^{2j}A(x)}{\partial x^{2j}}\right)\right].
\label{eq:diffusion_eq_sol_gen}
\end{align}
If we use $A(x)$ = sin($x$) for the initial condition in Eq. \eqref{eq:diffusion_eq_sol} we obtain
\begin{align}
\sum_{n=0}^{\infty}u^{n}(x,t) &= \sum_{n=0}^{\infty}\left(\frac{t^{n}}{n!}\right)\operatorname{exp}(-\beta t)\left[\sum^{n}_{j=0}\binom{n}{j}(\beta )^{n-j}c^{(j)}\left((-1)^{j}\operatorname{sin}(x)\right)\right].
\label{eq:diffusion_eq_sol_exp}
\end{align}
This method is trivially extendable to two-dimensional linear PDEs, the details of which are given in the Supplementary material (SM6). It is also possible to implement boundary conditions, and this is also demonstrated in the Supplementary material (SM7).

\subsection{Solving a nonlinear partial differential equation}

\begin{figure}[h!]
\centering
\begin{subfigure}[b]{0.4\textwidth}
	\includegraphics[width=1\textwidth]{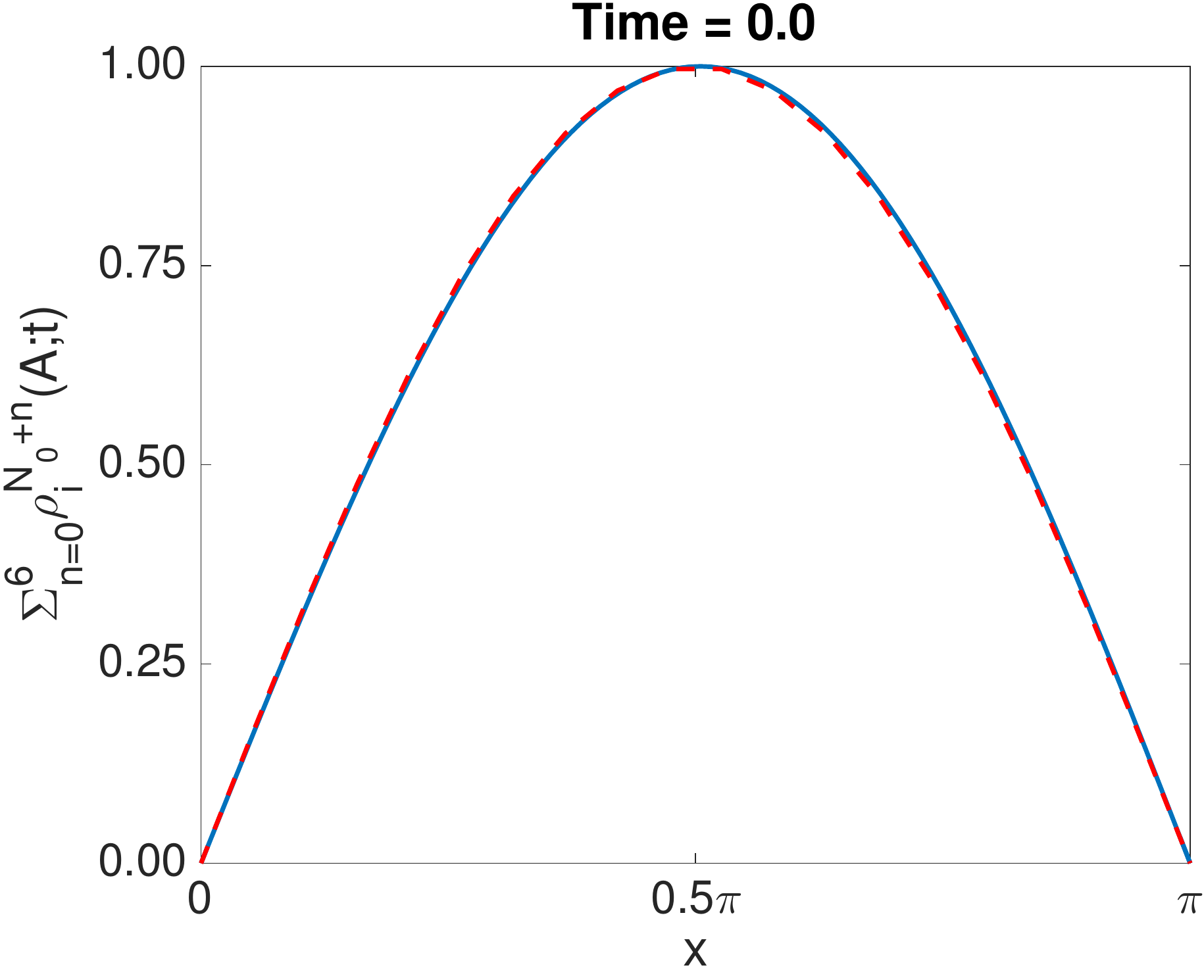}
	\subcaption{}
\end{subfigure}
\begin{subfigure}[b]{0.4\textwidth}
	\includegraphics[width=1\textwidth]{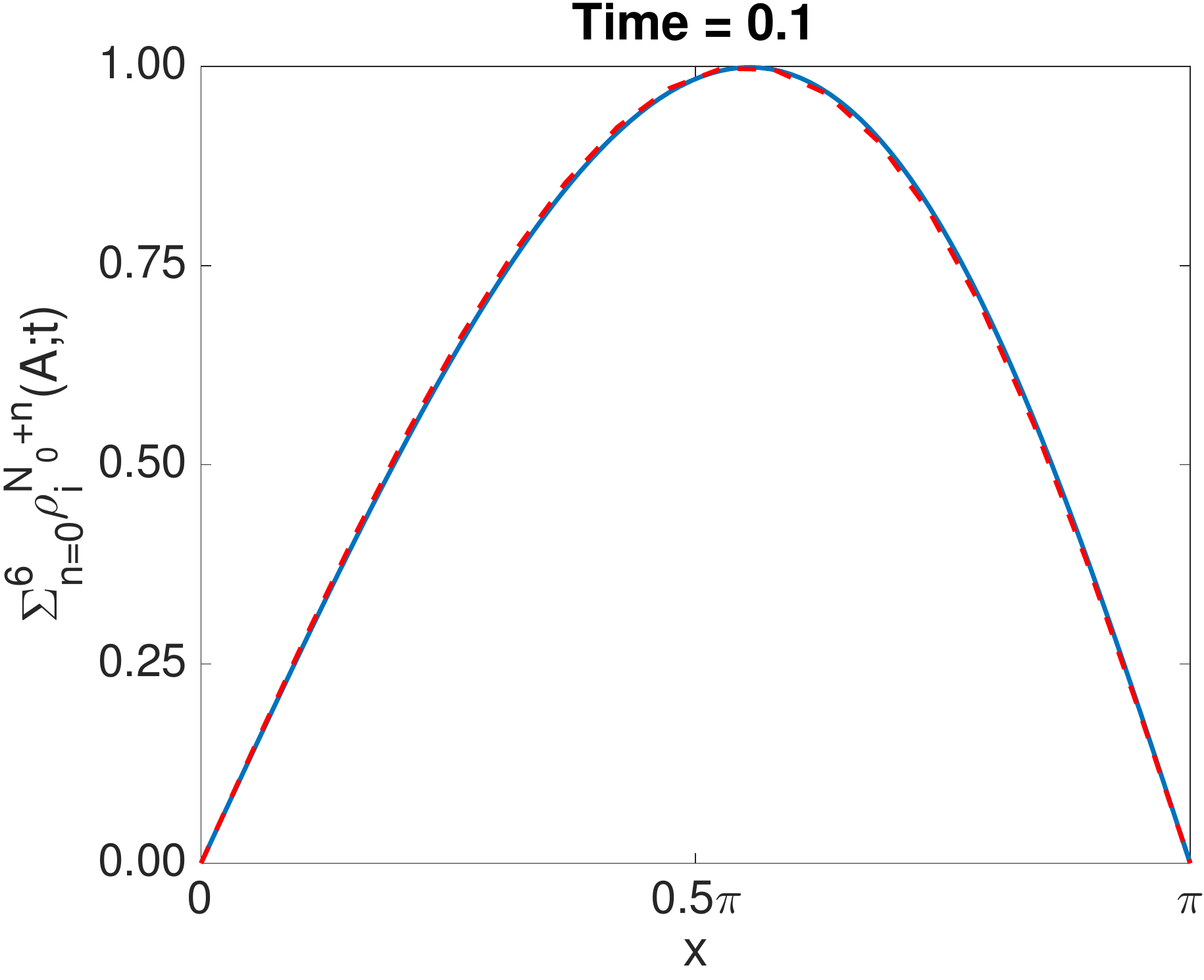}
	\subcaption{}
\end{subfigure}
\begin{subfigure}[b]{0.4\textwidth}
	\includegraphics[width=1\textwidth]{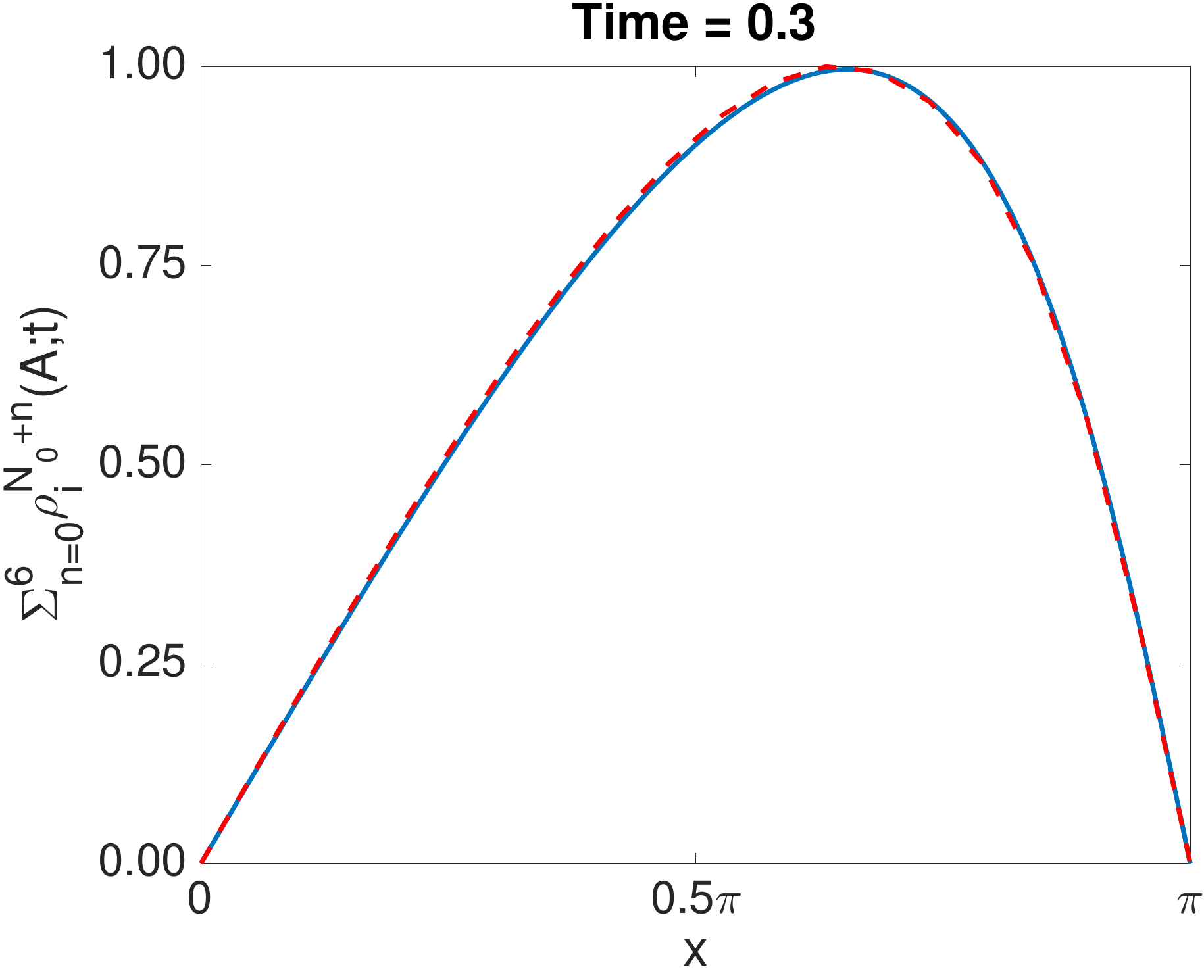}
	\subcaption{}
\end{subfigure}
\caption{A comparison of the solution to Eq. \eqref{eq:inviscid_burger_eq} and Eqs. \eqref{eq:inviscid_burger_eq_growth_simp_solve_L0}-\eqref{eq:inviscid_burger_eq_growth_simp_solve}. The blue lines indicate the solution of Eq. \eqref{eq:inviscid_burger_eq} and the red lines indicate the solutions of Eqs. \eqref{eq:inviscid_burger_eq_growth_simp_solve_L0}-\eqref{eq:inviscid_burger_eq_growth_simp_solve} for $\alpha = -0.5$. The initial condition $A(x)$ is $\operatorname{sin}(x) \in [0,\ \pi]$, and the truncation value for Eq. \eqref{eq:inviscid_burger_eq_growth_simp_solve} is $n=7$. In (a) $t = 0$, in (b) $t = 0.1$, and in (c) $t = 0.3$.}
\label{fig:figureNL}
\end{figure}

Finally, we demonstrate that this method is also extendable to nonlinear PDEs. For instance, the quasilinear inviscid Burgers equation.  The inviscid Burgers equation is
\begin{align}
\frac{\partial u}{\partial t} = \alpha u\frac{\partial u}{\partial x},
\label{eq:inviscid_burger_eq}
\end{align}
where $\alpha$ is a constant. We begin by writing Eq. \eqref{eq:inviscid_burger_eq} in the following manner
\begin{align}
\frac{\partial u^{n}}{\partial t} = -\beta u^{n} + \beta u^{n-1} + \alpha u^{n-1}\left(\sum^{\infty}_{i=0}\frac{\partial u^{i}}{\partial x}\right), \ \ \ \forall n > 0,
\label{eq:inviscid_burger_eq_growth}
\end{align}
with 
\begin{align}
\frac{\partial u^{0}}{\partial t} = -\beta u^{0}.
\label{eq:inviscid_burger_eq_growth0}
\end{align}
As in the case of nonlinear ODEs we have to multiply the $n^{th}$ stream by all other streams (including itself) to account for the nonlinearity in Eq. \eqref{eq:inviscid_burger_eq}.  We then decompose Eqs. \eqref{eq:inviscid_burger_eq_growth} and \eqref{eq:inviscid_burger_eq_growth0} in the following manner:
\begin{align}
\frac{\partial u^{0}}{\partial t} &= -\beta u^{0},
\label{eq:inviscid_burger_eq_growth_simp_solve_L0}
\end{align}
with
\begin{align}
\frac{\partial u^{1}}{\partial t} &= -\beta u^{1} + \beta u^{0} + \alpha u^{0}\frac{\partial u^{0}}{\partial x},
\label{eq:inviscid_burger_eq_growth_simp_solve_L1}
\end{align}
and
\begin{align}
\frac{\partial u^{n}}{\partial t} &= -\beta u^{n} + \beta u^{n-1} + \alpha u^{n-1}\left(\sum^{n-1}_{j=0}\frac{\partial u^{j}}{\partial x}\right) + \alpha\frac{\partial u^{n-1}}{\partial x}\left(\sum^{n-2}_{k=0}u^{k}\right), \ \ \ \forall n > 0.
\label{eq:inviscid_burger_eq_growth_simp_solve}
\end{align}
In Fig. \ref{fig:figureNL} the solution of Eqs. \eqref{eq:inviscid_burger_eq_growth_simp_solve_L0}-\eqref{eq:inviscid_burger_eq_growth_simp_solve} is compared with the solution of Eq. \eqref{eq:inviscid_burger_eq} before the onset of the multivalue behaviour that the solution of Eq. \eqref{eq:inviscid_burger_eq} exhibits. We use symbolic integration in Matlab to compute Eqs. \eqref{eq:inviscid_burger_eq_growth_simp_solve_L0}-\eqref{eq:inviscid_burger_eq_growth_simp_solve}.   It should be readily apparent how to extend this method to more complicated nonlinear PDEs.

\section{Discussion}

We have presented a power series method for solving both linear and nonlinear ODEs and PDEs. We finish by detailing some issues with the method we have introduced in this work.
\\
\\
Our main criticism of the work we have presented is that in the case of some nonlinear equations presented in this work we have not supplied solutions for the $n^{th}$ stream written in terms of simple functions. For instance, Eqs. \eqref{eq:non_lin2} and \eqref{eq:inviscid_burger_eq}. The method presented here would be most useful if an efficient means of writing the power series solutions for nonlinear equations became evident, which would allow analysis to be directly carried out on these solutions.  The Supplementary material (SM5) shows that in some cases of nonlinear equations it is possible to write the $n^{th}$ stream of its solution in terms of simple functions, however, a way to generalise the approach used on Eq. \eqref{eq:non_lin3} has not yet become apparent to the authors.
\\
\\
It is also important to acknowledge that we have not dealt with the issue of convergence in the power series we have presented.  It is obvious to say that the convergence of these power series, and their radius of convergence, will depend on the initial conditions of the equation, and the equation itself \citep{Abbott2001}.  However, a more general treatment on the convergence of the methods presented here is certainly required.  Finally, a word on the role of the shape parameter $\beta$.  Its role may seem somewhat superfluos, however, it is a simple way to circumvent numerical issues when the value of streams that compose solutions becomes too large for a standard computer to accurately represent. It also means that the value of the streams composing a solution can be made positive for a given interval of interest by selecting the appropriate value of $\beta$, and so provides another analytic tool to utilise when employing the methods presented here.

\section*{Acknowledgements}

RJHR would like to thank Kit Yates, Ruth Baker and Pierre Boutillier for helpful discussions. The author declares no competing interests.

\bibliography{References}

\newpage

\section*{Supplementary material}

\subsubsection*{SM1: The derivation of Equation \eqref{eq:static_diffusion} in the main text.}

We derive Eq. \eqref{eq:static_diffusion} in the following manner. The probability that an unbiased excluding random walker occupies site $i$ on a one-dimensional periodic lattice at time $t + \delta t$ is given by
\begin{align}
p_{i}(A;t+\delta t) &= p_{i}(A;t) + \frac{P_{m}\delta t}{2}\Bigg(p_{i-1,i}(A,0;t) - p_{i-1,i}(0,A;t)\Bigg) \nonumber \\
& \ \ \ + \frac{P_{m}\delta t}{2}\Bigg(p_{i,i+1}(0,A;t) - p_{i,i+1}(A,0;t)\Bigg).
\label{eq:SM11}
\end{align}
In Eq. \eqref{eq:SM11} $p_{i-1,i}(A,0;t)$ is the second-order probability that site $i-1$ and $i$ are occupied and unoccupied, respectively, at time $t$.  The other second-order terms in Eq. \eqref{eq:SM11} have similar meanings. If we rearrange Eq. \eqref{eq:SM11} and take $\delta t \rightarrow 0$ in the limit we obtain
\begin{align}
\frac{\mathrm{d}p_{i}(A;t)}{\mathrm{d}t} &= \frac{P_{m}}{2}\Bigg(p_{i-1,i}(A,0;t) - p_{i-1,i}(0,A;t)\Bigg) \nonumber \\
& \ \ \ + \frac{P_{m}}{2}\Bigg(p_{i,i+1}(0,A;t) - p_{i,i+1}(A,0;t)\Bigg).
\label{eq:SM12}
\end{align}
We now remove the second-order terms in Eq. \eqref{eq:SM12} by making the following closure
\begin{align}
p_{i,i+1}(A,0;t) = p_{i}(A;t)(1-p_{i+1}(A;t)).
\label{eq:SM13}
\end{align}
If we place Eq. \eqref{eq:SM13} in Eq. \eqref{eq:SM14} we obtain
\begin{align}
\frac{\mathrm{d}p_{i}(A;t)}{\mathrm{d}t} &= \frac{P_{m}}{2}\Bigg(p_{i-1}(A;t) - 2p_{i}(A;t) + p_{i+1}(A;t)\Bigg).
\label{eq:SM14}
\end{align}
If we drop the explicit `$A$' and `$t$' from our notation in Eq. \eqref{eq:SM14} we recapitulate Eq. \eqref{eq:static_diffusion}.

\subsubsection*{SM2: Algorithm for discrete random-walk}

We use a discrete random-walk model on a one-dimensional regular lattice with lattice spacing $\Delta$ \citep{Liggett} and length $N$, where $N$ is an integer describing the number of lattice sites.  Simulations are performed with either periodic boundary or no-flux conditions. Each random walker is assigned to a lattice site, from which it can move into an adjacent site. If an agent attempts to move into a site that is already occupied, the movement event is aborted. This process, whereby only one agent is allowed per site, is generally known as an exclusion process. Time is evolved continuously, and random walker movements are attempted in accordance with the Gillespie algorithm \citep{Gillespie_orig}. Attempted agent movement events occur with rate $P_{m}$ per unit time.  The initial conditions of the discrete model are provided in the main text when necessary.

\subsubsection*{SM3: The effect of different values of $\beta$ in Eq. \eqref{eq:metastatic_diffusion_sum} in the main text.}

\begin{figure}[h!]
\centering
\begin{subfigure}[b]{0.35\textwidth}
	\includegraphics[width=1\textwidth]{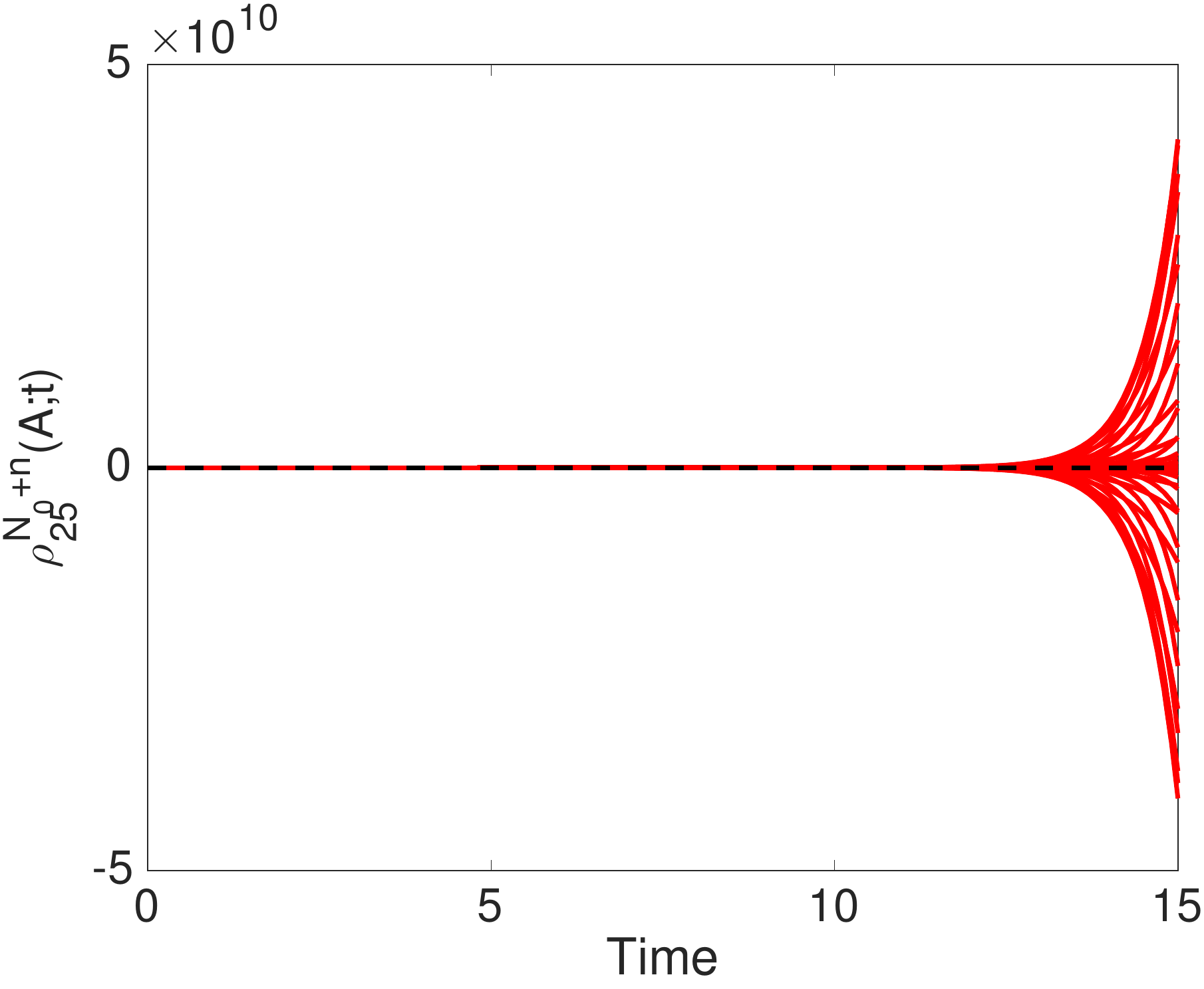}
	\subcaption{}
\end{subfigure}
\begin{subfigure}[b]{0.35\textwidth}
	\includegraphics[width=1\textwidth]{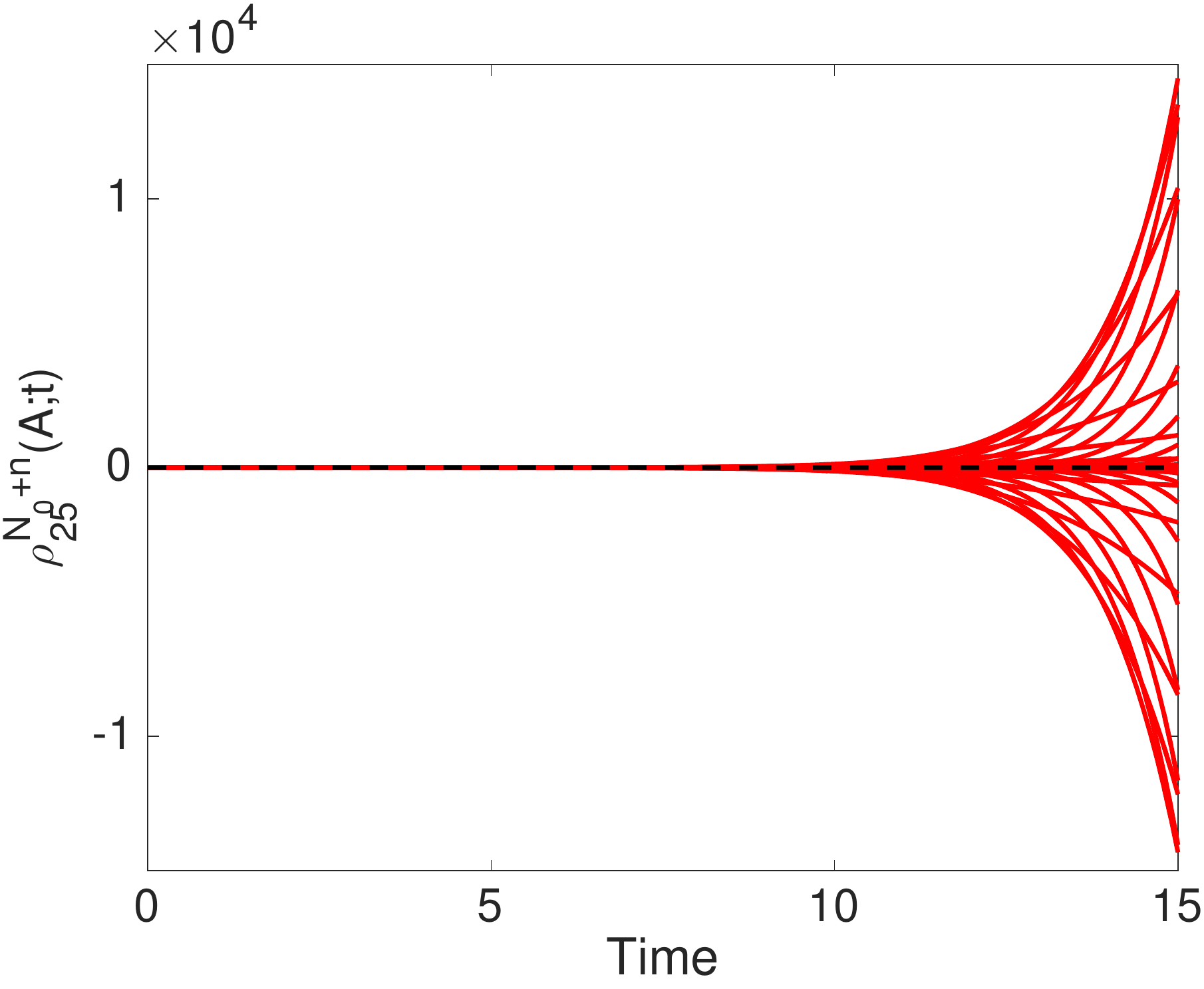}
	\subcaption{}
\end{subfigure}
\begin{subfigure}[b]{0.35\textwidth}
	\includegraphics[width=1\textwidth]{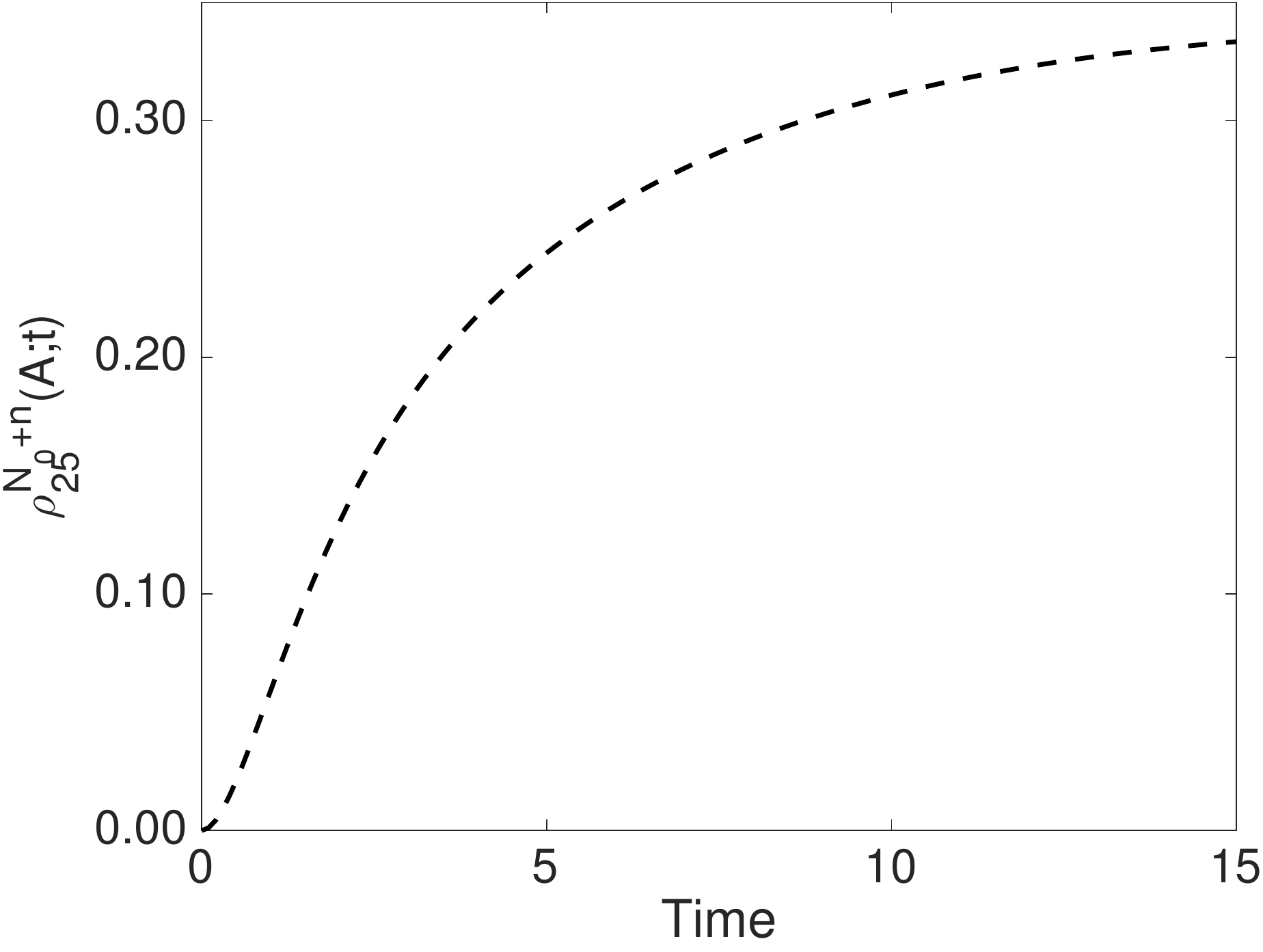}
	\subcaption{}
\end{subfigure}
\begin{subfigure}[b]{0.35\textwidth}
	\includegraphics[width=1\textwidth]{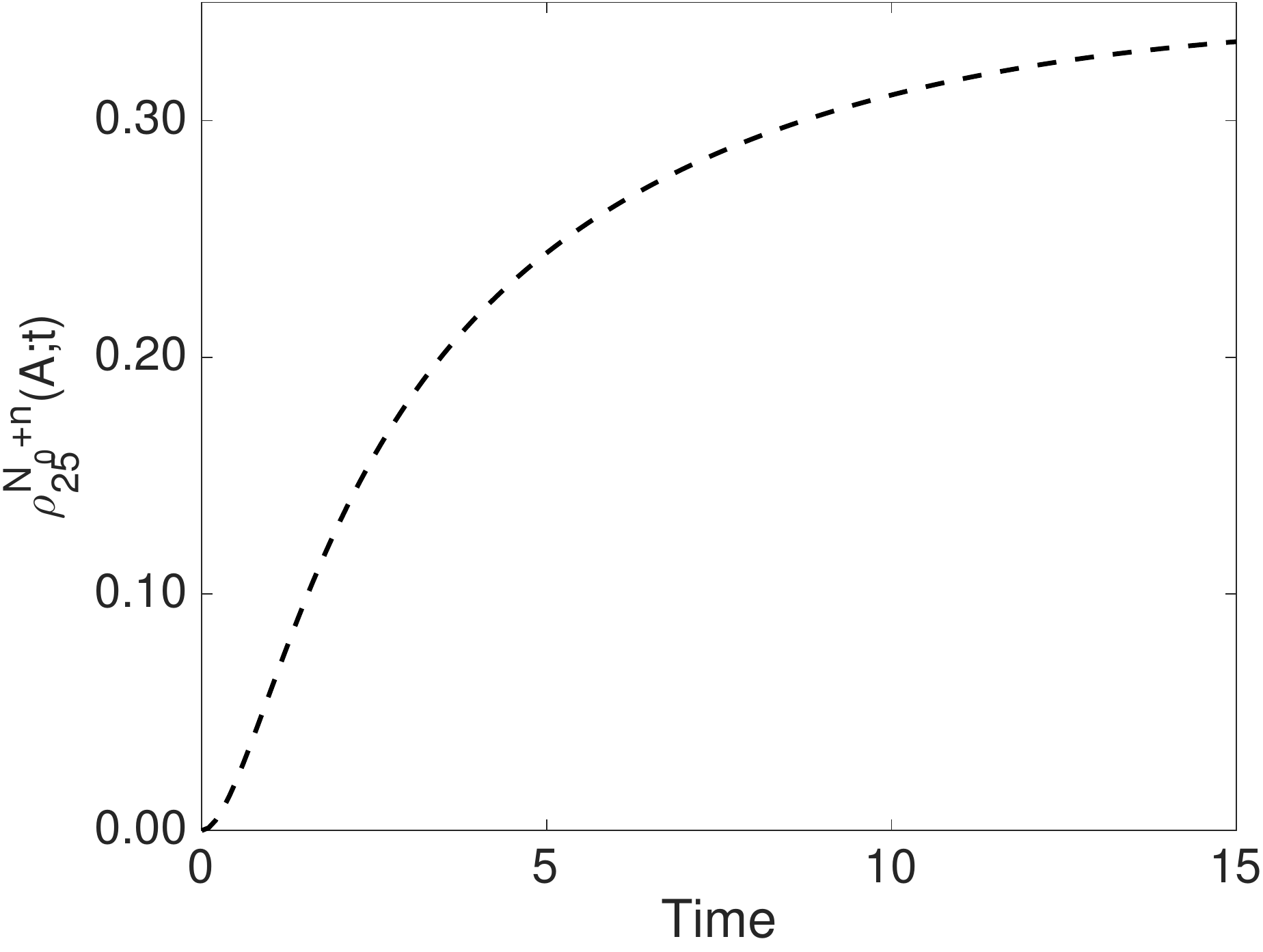}
	\subcaption{}
\end{subfigure}
\begin{subfigure}[b]{0.35\textwidth}
	\includegraphics[width=1\textwidth]{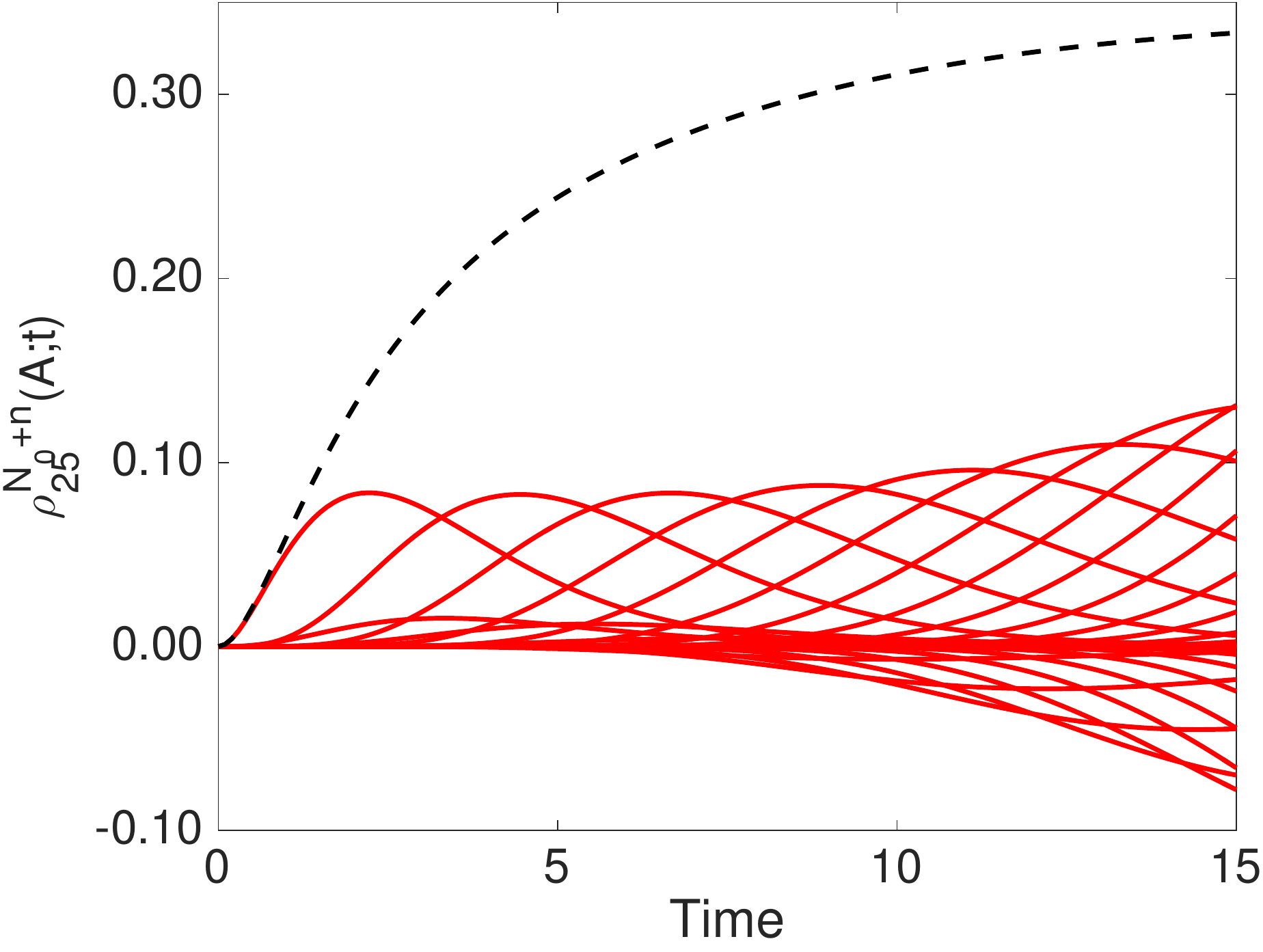}
	\subcaption{}
\end{subfigure}
\begin{subfigure}[b]{0.35\textwidth}
	\includegraphics[width=1\textwidth]{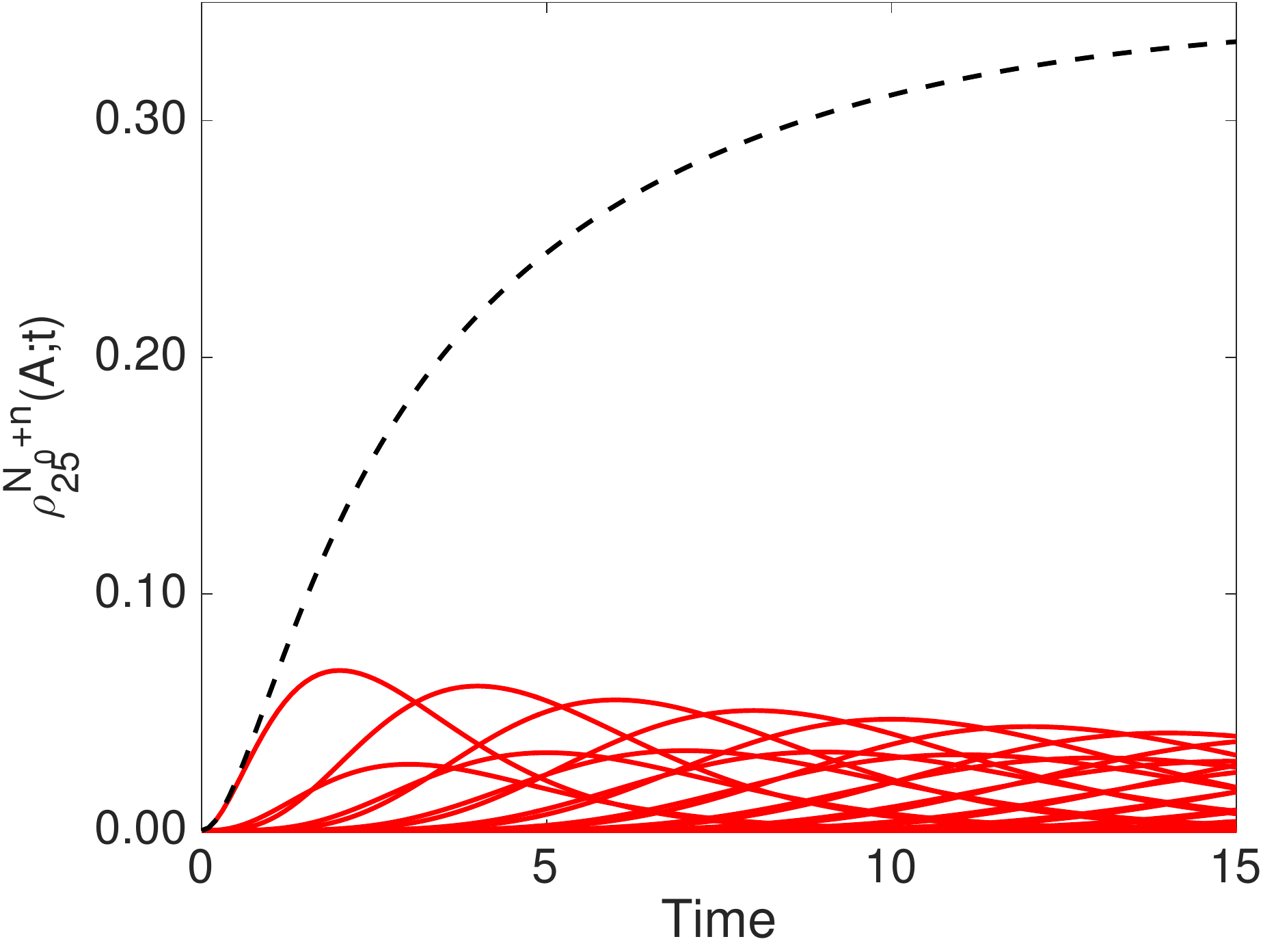}
	\subcaption{}
\end{subfigure}
\begin{subfigure}[b]{0.35\textwidth}
	\includegraphics[width=1\textwidth]{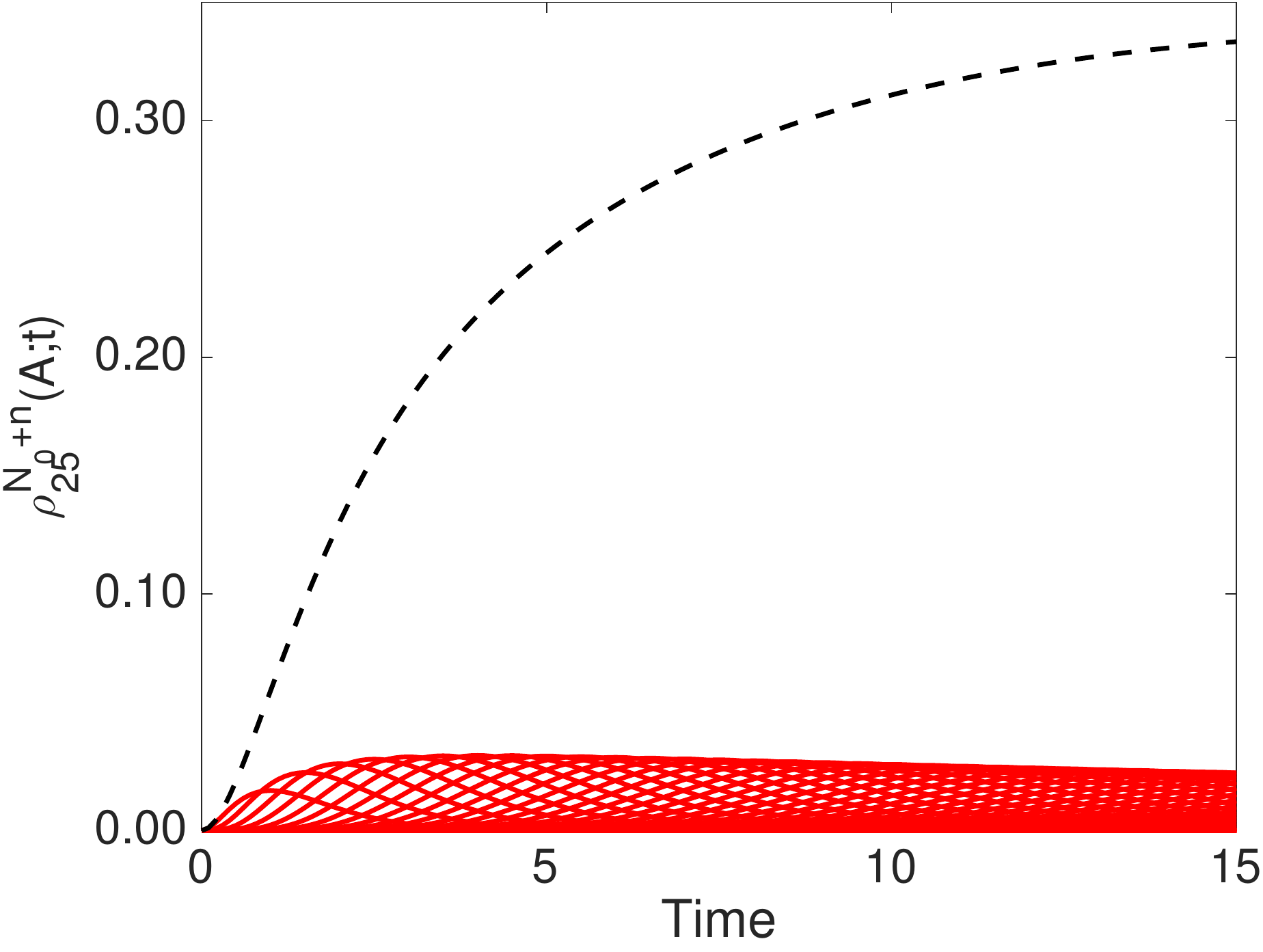}
	\subcaption{}
\end{subfigure}
\begin{subfigure}[b]{0.35\textwidth}
	\includegraphics[width=1\textwidth]{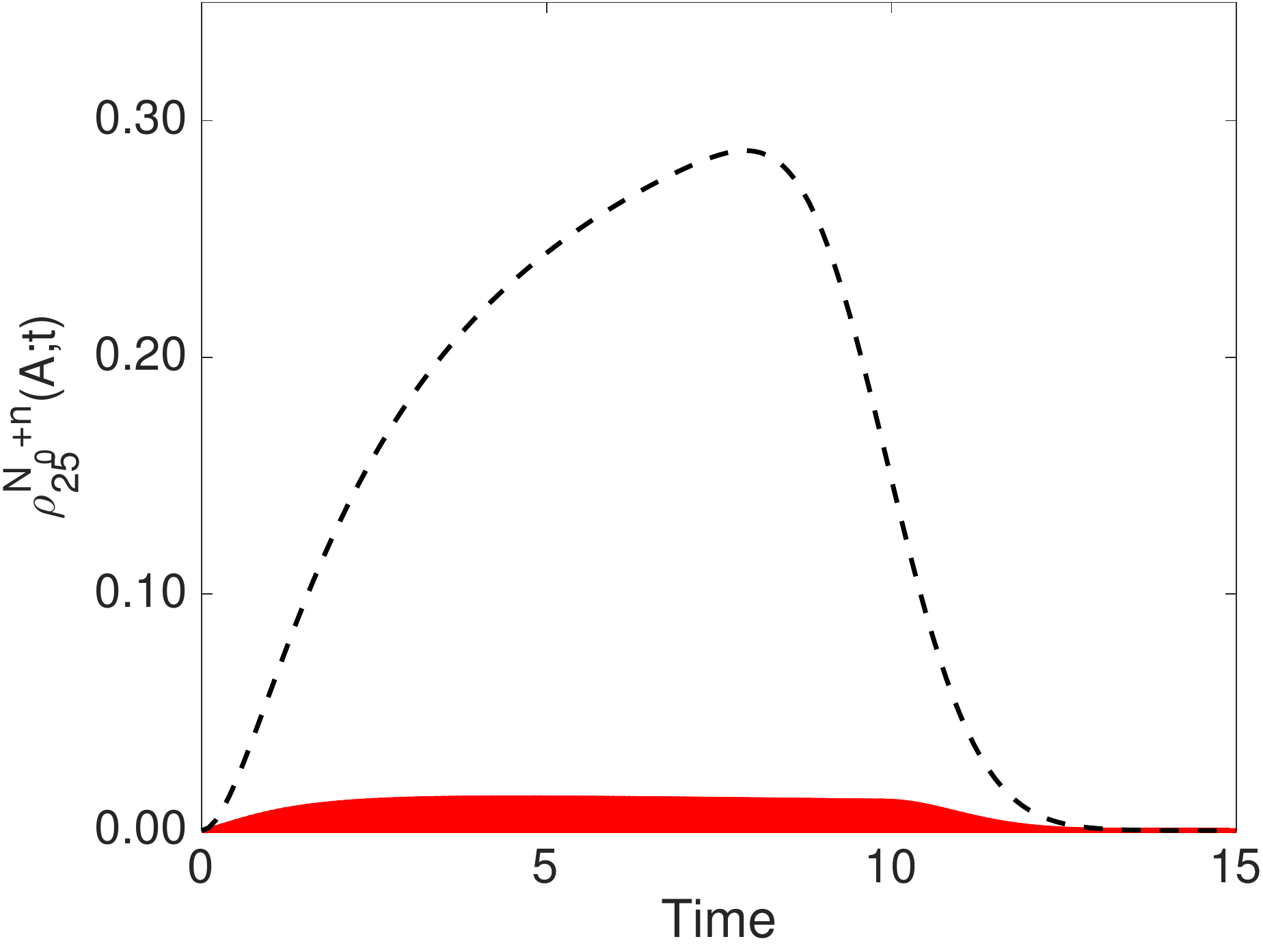}
	\subcaption{}
\end{subfigure}
\caption{The streams of site $i = 25$ as given by Eq. \eqref{eq:metastatic_diffusion_sol} for incrementing values of $n$ from 1:60 for different values of $\beta$. In (a) and (c) $\beta  = 0$, in (b) and (d) $\beta  = 0.5$, in (e) $\beta  = 0.9$, in (f) $\beta  = 1$, in (g) $\beta  = 2$, and in (h) $\beta  = 10$.  For all panels $P_{m} = 1$. It is evident that by selecting $\beta$ we are free to choose the shape of the streams. Panel (h) demonstrates what happens if the truncation of Eq. \eqref{eq:metastatic_diffusion_sum} is too low (the truncation value in this case is 60), and/or $\beta$ is too large. The black-dashed line is the sum of all the streams for site $i=25$, given by Eq. \eqref{eq:metastatic_diffusion_sum}. The solution given by Eq. \eqref{eq:metastatic_diffusion_sum} is the same for panels (a)-(g), but fails in panel (h) for the reasons discussed.}
\label{fig:figure2c}
\end{figure}

In Fig. \ref{fig:figure2c} we display streams, Eq. \eqref{eq:metastatic_diffusion_sol}, for different values of $\beta$. This demonstrates how $\beta $ influences the shape of the streams that compose the solution given by Eq. \eqref{eq:metastatic_diffusion_sum}.

\subsubsection*{SM4: Linear ordinary differential equation}

We provide the solution to the following linear ODE
\begin{align}
\frac{\mathrm{d}q}{\mathrm{d}t} = (1-2t)q,
\label{eq:non_lin1}
\end{align}
which is linear in $q$. The analytic solution of Eq. \eqref{eq:non_lin1} is
\begin{align}
q(t) = C_{2}\operatorname{exp}(t-t^{2}),
\label{eq:non_lin1_sol}
\end{align}
where $C_{2}$ is the value of $q(t)$ at $t = 0$.
To solve Eq. \eqref{eq:non_lin1} in our framework we rewrite it as
\begin{align}
\frac{\mathrm{d}q^{n}(t)}{\mathrm{d}t} = -\beta q^{n}(t) + \beta q^{n-1}(t) + (1-2t)q^{n-1}(t), \ \ \ \forall \ n > 0,
\label{eq:nonlin_exp}
\end{align}
with 
\begin{align}
\frac{\mathrm{d}q^{0}(t)}{\mathrm{d}t} = -\beta q^{0}(t).
\label{eq:nonlin_exp0}
\end{align}
In Fig. \ref{fig:figureNL} (a) the solution of Eqs. \eqref{eq:nonlin_exp} and \eqref{eq:nonlin_exp0} is compared with the analytical solution Eq. \eqref{eq:non_lin1_sol}.

\begin{figure}[h!]
\centering
\begin{subfigure}[b]{0.4\textwidth}
	\includegraphics[width=1\textwidth]{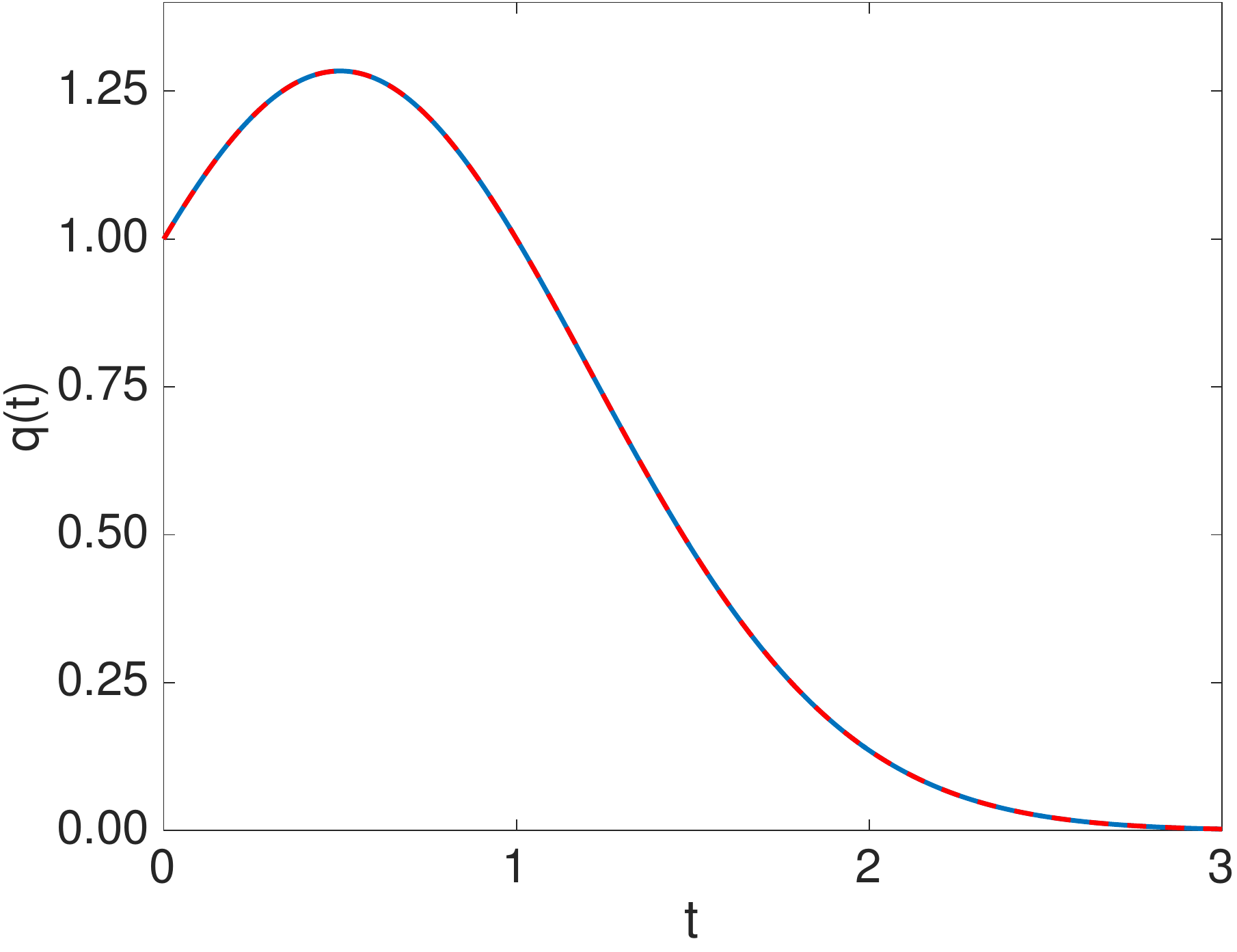}
	\subcaption{}
\end{subfigure}
\caption{In (a) Eqs. \eqref{eq:nonlin_exp} and \eqref{eq:nonlin_exp0} are compared with Eq. \eqref{eq:non_lin1_sol} for $C_{2} = 1$. The truncation value for Eq. \eqref{eq:nonlin_exp} is $n=30$, and $\beta = 10$.}
\label{fig:figureLMSM}
\end{figure}

\subsubsection*{SM5: Nonlinear ordinary differential equation}

To solve Eq. \eqref{eq:non_lin3} we proceed in the following manner. We begin with
\begin{align}
\frac{\mathrm{d}y}{\mathrm{d}t} = \alpha y^{2},
\label{eq:SM31}
\end{align}
and
\begin{align}
y^{0} = A,
\label{eq:SM32}
\end{align}
and
\begin{align}
\frac{\mathrm{d}y^{1}}{\mathrm{d}t} = \alpha yy^{0}.
\label{eq:SM33}
\end{align}
In this derivation we assume $\beta = 0$ for simplicity.
Initially, we rewrite Eq. \eqref{eq:SM31} as
\begin{align}
\frac{\mathrm{d}}{\mathrm{d}t}(\operatorname{log}(y)) = \alpha y.
\label{eq:SM34}
\end{align}
This means
\begin{align}
\frac{\mathrm{d}y^{1}}{\mathrm{d}t} = y^{0}\frac{\mathrm{d}}{\mathrm{d}t}(\operatorname{log}(y)),
\label{eq:SM35}
\end{align}
which gives
\begin{align}
y^{1} = y^{0}\operatorname{log}(y) + c_{1}.
\label{eq:SM36}
\end{align}
Therefore
\begin{align}
y = y^{0}\operatorname{exp}\left(\frac{y^{1}}{y^{0}}\right),
\label{eq:SM37}
\end{align}
because $c_{1} = -y^{0}\operatorname{log}(y^{0})$. If we place Eq. \eqref{eq:SM37} in Eq. \eqref{eq:SM33} we obtain
\begin{align}
\frac{\mathrm{d}y^{1}}{\mathrm{d}t} = \alpha (y^{0})^{2}\operatorname{exp}\left(\frac{y^{1}}{y^{0}}\right),
\label{eq:SM38}
\end{align}
which we can integrate to obtain
\begin{align}
y^{1} = -y^{0}\operatorname{log}\left(1 - \alpha y^{0}t\right).
\label{eq:SM39}
\end{align}
Now we recognise 
\begin{align}
\frac{\mathrm{d}y^{2}}{\mathrm{d}t} = \frac{y^{1}}{y^{0}}\frac{\mathrm{d}y^{1}}{\mathrm{d}t}.
\label{eq:SM310}
\end{align}
If we integrate Eq. \eqref{eq:SM310}, and then solve for $y^{3}$ in a similar manner we obtain the following power series solution for $y$
\begin{align}
y = \sum^{\infty}_{n=0} \frac{(-1)^{n}y^{0}}{n!}\operatorname{log}^{n}(1 - \alpha y^{0} t).
\label{eq:non_lin3solSM}
\end{align}
Alternatively, we can place Eq. \eqref{eq:SM39} into Eq. \eqref{eq:SM37} to obtain
\begin{align}
y = \frac{y^{0}}{1 - \alpha y^{0} t},
\label{eq:non_lin3analSM}
\end{align}
which recapitulates the analytical solution to Eq. \eqref{eq:SM31}, given as Eq. \eqref{eq:non_lin3anal}. 

\subsubsection*{SM6: Two-dimensional linear partial differential equation}

The general solution for the two-dimensional linear diffusion equation $(\beta = 0)$ is
\begin{align}
\sum^{\infty}_{n=0}p^{n}(x,y;t) = \sum^{\infty}_{n=0}\frac{(Dt)^{n}}{n!}\left[\sum^{n}_{i}\binom{n}{i}\frac{\partial^{2n}A(x,y)}{\partial x^{2(n-i)}\partial y^{2i}} \right].
\label{eq:lin_bc}
\end{align}

\subsubsection*{SM7: Boundary conditions for linear partial differential equation}

We now demonstrate how to implement boundary conditions in linear PDEs with our method. A simple example is for the diffusion equation
\begin{align}
\frac{\partial u}{\partial t} = D\frac{\partial^{2} u}{\partial x^{2}},
\label{eq:lin_bc}
\end{align}
with
\begin{align}
u^{L}(x,0) = A(x) = \gamma + \lambda x,
\end{align}
and
\begin{align}
\sum^{\infty}_{i=0}{u^{L+i \delta L}(0,t)} = \gamma = u^{L+i \delta L}(0,0), \ \ \ \ \sum^{\infty}_{i=0}{u^{L+i \delta L}(L,t)} = \gamma + \lambda L = u^{L+i \delta L}(L,0).
\end{align}
If we implemented Neumann boundary conditions these would take the form
\begin{align}
\sum^{\infty}_{i=0}\frac{\partial u^{L+i \delta L}(0,t)}{\partial x} = -c_{1}, \ \ \ \ \sum^{\infty}_{i=0}\frac{\partial u^{L+i \delta L}(L,t)}{\partial x} = c_{2}.
\end{align}
To implement our boundary conditions we proceed as follows: As we already have the general solution for Eq. \eqref{eq:diffusion_eq_sol} we can take its partial derivative with respect to time to obtain\footnote{A simpler way to obtain the solution for the given initial condition is the following. Implementing initial condition in Eq. \eqref{eq:lin_bc} gives:
\begin{align}
\frac{\partial u^{L+\delta L}}{\partial t} = -\beta u^{L + \delta L} + \beta u^{L},
\end{align} 
with the (straightforward) solution being
\begin{align}
\sum^{\infty}_{i=0}u^{L+i\delta L}(x,t) = \sum^{\infty}_{i=0}A(x)\frac{(\beta t)^{i}}{i!}\operatorname{exp}(-\beta t) = A(x)\sum^{\infty}_{i=0}\frac{(\beta t)^{i}}{i!}\operatorname{exp}(-\beta t) = A(x),
\end{align} 
as the Poisson distribution sums to identity.}
\begin{align}
\frac{\partial u}{\partial t} = \sum_{n=0}^{\infty}\left(n\left(\frac{t^{n-1}}{n!}\right)\operatorname{exp}(-\beta t)-\beta \left(\frac{t^{n}}{n!}\right)\operatorname{exp}(-\beta t)\right)\left[\sum^{n}_{j=0}\binom{n}{j}(\beta )^{n-j}D^{(j)}\left(\frac{\partial^{2j}A(x)}{\partial x^{2j}}\right)\right],
\end{align}
which means
\begin{align}
\frac{\partial^{2} u}{\partial x^{2}} &= \nonumber \\
& \hspace{-0cm} \left(\frac{1}{D}\right)\sum_{n=0}^{\infty}\left(n\left(\frac{t^{n-1}}{n!}\right)\operatorname{exp}(-\beta t)-\beta \left(\frac{t^{n}}{n!}\right)\operatorname{exp}(-\beta t)\right)\left[\sum^{n}_{j=0}\binom{n}{j}(\beta )^{n-j}D^{(j)}\left(\frac{\partial^{2j}A(x)}{\partial x^{2j}}\right)\right].
\label{eq:diffusion_dx2}
\end{align}
Integrating Eq. \eqref{eq:diffusion_dx2} with respect to $x$ gives
\begin{align}
\frac{\partial u}{\partial x} + c_{1} &= \nonumber \\
&\hspace{-1cm} \left(\frac{1}{D}\right)\sum_{n=0}^{\infty}\left(n\left(\frac{t^{n-1}}{n!}\right)\operatorname{exp}(-\beta t)-\beta \left(\frac{t^{n}}{n!}\right)\operatorname{exp}(-\beta t)\right)\left[\sum^{n}_{j=0}\binom{n}{j}(\beta )^{n-j}D^{(j)}\left(\frac{\partial^{2j-1}A(x)}{\partial x^{2j-1}}\right)\right],
\label{eq:diffusion_dx}
\end{align}
and integrating Eq. \eqref{eq:diffusion_dx} with respect to $x$ gives
\begin{align}
u + c_{1}x + c_{2} &= \nonumber \\
&\hspace{-2cm} \left(\frac{1}{D}\right)\sum_{n=0}^{\infty}\left(n\left(\frac{t^{n-1}}{n!}\right)\operatorname{exp}(-\beta t)-\beta \left(\frac{t^{n}}{n!}\right)\operatorname{exp}(-\beta t)\right)\left[\sum^{n}_{j=0}\binom{n}{j}(\beta )^{n-j}D^{(j)}\left(\frac{\partial^{2j-2}A(x)}{\partial x^{2j-2}}\right)\right].
\label{eq:full_sol}
\end{align}
If we apply the boundary conditions and initial condition to Eq. \eqref{eq:full_sol} we obtain
\begin{align}
u &= \frac{1}{D}\sum_{n=0}^{\infty}\operatorname{exp}(-\beta t)\left(n\left(\frac{t^{n-1}}{n!}\right)-\beta \left(\frac{t^{n}}{n!}\right)\right)\Bigg[\left[(\beta )^{n}\left(\gamma \frac{x^{2}}{2} + \lambda \frac{x^{3}}{6}\right) + n(\beta )^{n-1}D(\gamma + \lambda x)\right] \nonumber \\
& \ \ \ - \left(\frac{x}{L}\right)\left[(\beta )^{n}\left(\gamma \frac{L^{2}}{2} + \lambda \frac{L^{3}}{6}\right) + n(\beta )^{n-1}D(\gamma + \lambda L)\right] + \left(\frac{x}{L}-1\right)\left[n(\beta )^{n-1}D\gamma\right]\Bigg]  \nonumber \\ & \ \ \ + \lambda x + \gamma
\label{eq:lin_bc_sol}
\end{align}
as the solution to Eq. \eqref{eq:lin_bc}. In this solution we define $$n\left(\frac{t^{n-1}}{n!}\right) = 0$$ when $n = 0$ to avoid a singularity when $t = 0$.

\end{document}